\documentclass[final,twoside]{siamart1116}

\numberwithin{theorem}{section}
\numberwithin{equation}{section}

\title{
  Generalized Sinkhorn iterations for regularizing inverse problems using optimal mass transport%
  \thanks{This work was supported by the Swedish Research Council (VR) grant 2014-5870, the Swedish Foundation of Strategic Research (SSF) grant AM13-0049, National Science Foundation (NSF) grant CCF-1218388, and the Center for Industrial and Applied Mathematics (CIAM).}
}

\author{
  Johan Karlsson%
  \thanks{Division of Optimization and Systems Theory, Department of Mathematics, KTH  Royal Institute of Technology, 100 44 Stockholm, Sweden. Contact: \email{johan.karlsson@math.kth.se}, \email{aringh@kth.se}.}
  \and
  Axel Ringh%
  \footnotemark[2]
}

\newenvironment{example}{\refstepcounter{theorem} {\em Example} \thetheorem.}{}
\newenvironment{remark}{\refstepcounter{theorem} {\em Remark} \thetheorem.}{}
\newenvironment{assumption}{\refstepcounter{theorem} {\em Assumption} \thetheorem.}{}

\usepackage{amsmath}
\usepackage{amssymb}
\usepackage{amstext}
\usepackage{graphicx}

\usepackage{algorithm}
\usepackage{algorithmic}

\usepackage{color}

\usepackage{url}

\usepackage{psfrag}

\usepackage[caption=false]{subfig}

\usepackage{epstopdf}

\usepackage{booktabs}

\usepackage{dsfont}

\usepackage{bbm}

\usepackage{bm}

\usepackage{yfonts}

\usepackage{wasysym}

\usepackage{euscript}

\usepackage{latexsym}

\newcommand{\tr}{{\rm trace}}
\newcommand{\diag}{{\rm diag}}
\newcommand{\mR}{{\mathbb R}}
\newcommand{\cI}{{\mathcal I}}
\newcommand{\mue}{{\mu_{\rm est}}}
\newcommand{\mut}{{\mu_{\rm true}}}
\newcommand{\ett}{{\bf 1}}
\DeclareMathOperator*{\argmin}{arg \, min}
\DeclareMathOperator{\Prox}{Prox}

\def \myFigPathSheppLogan {2D_tomo_shepp_logan/}

\def \SheppLoganPrior {"Prior_no_title"}
\def \SheppLoganPhantom {"Phantom_no_title"}
\def \SheppLoganFBP {"Filtered_Backprojection_no_title"}
\def \SheppLoganTV {"TV_reconstruction_no_title"}
\def \SheppLoganL {"L2_plus_TV_regularization_reg_param_10_0_no_title"}
\def \SheppLoganLa {"L2_plus_TV_regularization_reg_param_1_0_no_title"}
\def \SheppLoganLb {"L2_plus_TV_regularization_reg_param_100_0_no_title"}
\def \SheppLoganLc {"L2_plus_TV_regularization_reg_param_10000_0_no_title"}
\def \SheppLoganOMT {"Optimal_transport_TV_reconstruction_reg_param_4_0_no_title"}

\def \myFigPathHand {2D_tomo_hand/}

\def \HandPrior {"Prior_no_title"}
\def \HandPhantom {"Phantom_no_title"}
\def \HandFBP {"Filtered_Backprojection_no_title"}
\def \HandTV {"TV_reconstruction_no_title"}
\def \HandLb {"L2_plus_TV_regularization_reg_param_10_0_no_title"}
\def \HandOMT {"Optimal_transport_TV_reconstruction_reg_param_4_0_no_title"}
\def \HandMaskedPrior {"masked_prior_no_title"}
\def \HandOMTMasked {"Optimal_transport_TV_reconstruction_reg_param_4_0mass_movment_postManipulation_no_title"}

\begin{document}

\maketitle

\begin{abstract}
The optimal mass transport problem gives a geometric framework for optimal allocation, and has recently gained significant interest in application areas such as signal processing, image processing, and computer vision. 
Even though it can be formulated as a linear programming problem, it is in many cases intractable for large problems due to the vast number of variables. 
A recent development to address this builds on an approximation with an entropic barrier term and solves the resulting optimization problem using Sinkhorn iterations. 
In this work we extend this methodology to a class of inverse problems. 
In particular we show that Sinkhorn-type iterations can be used to compute the proximal operator of the transport problem for large problems. 
A splitting framework is then used to solve inverse problems where the optimal mass transport cost is used for incorporating a priori information.  
We illustrate the method on problems in computerized tomography. In particular we consider a  limited-angle computerized tomography problem, where a priori information is used to compensate for missing measurements.

\end{abstract}

\section{Introduction}

The optimal mass transport problem provides a useful framework that can be utilized in many contexts, ranging from optimal allocation of resources to applications in imaging and machine learning. 
In this work we extend this framework and consider optimization problems where the objective function contains an optimal mass transport cost.
This includes several inverse problems of interest, for example to model deformations in the underlying object. 

The optimal mass transport problem is sometimes referred to as the Monge-Kantorovich transportation problem 
and 
was originally formulated by Gaspard Monge in 1781 \cite{villani2008optimal} 
for transport of soil for construction of forts and roads in order to minimize the transport expenses. He described the problem as follows: ``divide two equal volumes into infinitesimal particles and associate them one to another so that the sum of the path lengths multiplied by the volumes of the particles be minimum possible'' \cite{kantorovich2006problem}. 
The second founder of the field is the mathematician and economist Leonid Kantorovich, who made major advances to the area and as part of this reformulated the problem as a convex optimization problem along with a dual framework.
He later received the Nobel Memorial Prize in Economic Sciences for his contributions to the theory of optimum allocation of resources \cite{vershik2013long}.
For an introduction and an overview of the optimal mass transport problem, see, e.g., \cite{villani2008optimal}.
During the last few decades the approach has gained much interest in several application fields such as image processing \cite{haker2004optimal, de2011optimal,kaijser1998computing}, signal processing \cite{engquist2014application, georgiou2009metrics, jiang2012geometric}, computer vision and machine learning \cite{rubner2000earth, ling2007efficient, carli2013convex,sadeghian2015automatic}. 

In our setting the optimal transportation cost is used as a distance for comparing objects and incorporating a priori information. An important property of this distance is that it does not only compare objects point by point, as standard $L^p$ metrics, but instead quantifies the length with which that mass is moved. This property makes the distance natural for quantifying uncertainty and modelling deformations \cite{jiang2012geometric,karlsson2013uncertainty,karlsson2014robustness}.
More specifically geodesics (in, e.g., the associated Wasserstein-2 metric \cite{villani2008optimal}) preserve ``lumpiness,'' and when linking objects via
geodesics of the metric there is a natural deformation between the objects.
Such a property appears highly desirable in tracking moving objects and integrating data
from a variety of sources (see, e.g., \cite{jiang2012geometric}).
This is in contrast to the fade-in-fade-out effect that occurs when linking objects in standard metrics (e.g., the $L^2$ metric). 

Although the optimal mass transport problem has many desirable properties it also has drawbacks. 
Monge's formulation is a nonconvex optimization problem and the Kantorovich formulation results in large-scale optimization problems that are intractable to solve with standard methods even for modest size transportation problems.
A recent development to address this computational problem builds on adding an entropic barrier term and solving the resulting optimization problem using the so called Sinkhorn iterations \cite{cuturi2013sinkhorn}. This allows for computing an approximate solution of large transportation problems and
has proved useful for many problems where no computationally feasible method previously exists. Examples include computing multi-marginal optimal transport problems and barycenters (centroids) \cite{benamou2015iterative} and sampling from multivariate probability distribution \cite{jacob2016coupling}.  For quadratic cost functions this approach can also be seen as the solution to a Schr\"odinger bridge problem \cite{chen2016relation}. 

In this work we build on these methods and consider variational problems that contain an optimal mass transport cost. In particular, we focus on problems of the form 
\begin{align}\label{eq:Pgen0}
\min_{\mue} &\quad  T_\epsilon(\mu_0, \mue)+g(\mue),
\end{align} 
where $T_{\epsilon}$ is the entropy regularized optimal mass transport cost, and $g$ both quantify data miss-match and can contain other regularization terms. 
Typically $\mu_0$ is a given prior and the minimizing argument $\mue$ is the sought reconstruction.
This formulation allows us to model deformations, to address problems with unbalanced masses (cf. \cite{georgiou2009metrics}), as well as solving  gradient flow problems \cite{benamou2016anaugmented, peyre2015entropic}
appearing in the Jordan-Kinderlehrer-Otto framework \cite{jordan1998variational}.
A common technique for solving optimization problems with several additive terms is to utilize the proximal point method \cite{rockafellar1976monotone} together with a variable splitting technique \cite{bauschke2011convex, boyd2011distributed,eckstein1989splitting,komodakis2015playing}. These methods typically utilize the proximal operator and provide a powerful computational framework for solving composite optimization problems. 
Similar frameworks have previously been used for solving optimization problems that include a transportation cost. In \cite{benamou2016anaugmented} a fluid dynamics formulation   
\cite{benamou2000computational} is considered that can be used to compute the proximal operator, and in \cite{peyre2015entropic} an entropic proximal operator \cite{teboulle1992entropic} is used for solving problems of the form \eqref{eq:Pgen0}.
In this paper we propose a new fast iterative computational method for computing the proximal operator of $T_\epsilon(\mu_0,\cdot)$ based on Sinkhorn iterations.
This allows us to use splitting methods, such as Douglas-Rachford type methods \cite{eckstein1992douglas,bot2013douglas}, in order to solve large scale inverse problems of interest in medical imaging.

The paper is structured as follows. 
In Section~\ref{sec:background} we describe relevant background material on optimal transport and Sinkhorn iterations, and we also review the proximal point algorithm and variable splitting in optimization.
Section~\ref{sec:dual} considers the dual problem of \eqref{eq:Pgen0} and we introduce generalized Sinkhorn iterations which can be used for computing \eqref{eq:Pgen0} efficiently for several cost functions $g$. In particular the proximal operator of $T_\epsilon(\mu_0,\cdot)$ is computed using iterations with the same computational cost as Sinkhorn iterations.
In Section~\ref{sec:inverse} we describe how a Douglas-Rachford type splitting technique can be applied for solving a general class of large scale problems on the form \eqref{eq:Pgen0} and in Section~\ref{sec:applications} we apply these techniques to inverse problems using an optimal transport prior, namely two reconstruction problems in Computerized Tomography. In Section~\ref{sec:conclusions} we discuss conclusions and further directions. Finally, the paper has four appendices containing deferred proofs, connections with the previous work \cite{peyre2015entropic}, as well as details regarding the numerical simulations.

\subsection{Notation}
Next we briefly introduce some notation. Most operations in this paper are defined elementwise. In particular we use $\odot$, $./$, $\exp$, and $\log$, to denote elementwise multiplication, division, exponential function, and logarithm function, respectively. We also use $\le$ ($<$) to denote elementwise  inequality (strict).
Let $\cI_S(\mu)$ be the indicator function of the set $S$, i.e., $\cI_S(\mu)=0$ if $\mu\in S$ and $\cI_S(\mu)=+\infty$ otherwise.
Finally, let $\ett_{n}$ denote the $n\times 1$ (column) vector of ones. 

\section{Background}\label{sec:background}

\subsection{The optimal mass transport problem and entropy regularization}

The Monge's formulation of the optimal mass transport problem is set up as follows. Given two nonnegative functions, $\mu_0$ and $\mu_1$, of the same mass, defined on a compact set $X\subset \mR^d$, one seeks the transport function $\phi: X\to X$ 
that minimizes the transportation cost
\[
\int_{X}c(\phi(x),x)\mu_0(x)dx
\]
and that is a mass preserving map from $\mu_0$ to $\mu_1$, i.e., 
\[
\int_{x\in A}\mu_{1}(x)dx=\int_{\phi(x)\in A}\mu_{0}(x)dx\quad \mbox{ for all } A \subset X.
\]
Here $c(x_0,x_1):X\times X\to \mR_+$ is a cost function that describes the cost for transporting a unit mass from $x_0$ to $x_1$. It should be noted that this optimization problem is not convex and the formulation is not symmetric with respect to functions $\mu_0$ and $\mu_1$.

In the Kantorovich formulation one instead introduce a transference plan, $M:X\times X\to \mR_+$, which characterizes the mass which is moved from $x_0$ to $x_1$. This construction generalizes to general nonnegative measures $dM\in {\mathcal M}_+(X\times X)$ and allows for transference plans where the mass in one point in $\mu_0$ is transported to a set of points in $\mu_1$. The resulting optimization problem is convex and the cost is given by
\begin{align}
T(\mu_0, \mu_1) =
\min_{dM\in {\mathcal M}_+(X\times X)}&
\int_{(x_0,x_1)\in X\times X}c(x_0,x_1)dM(x_0,x_1)\nonumber\\
\mbox{subject to }\; &\mu_0(x_0)dx_0=\int_{x_1\in X}dM(x_0,x_1),\label{eq:Kantorovich}\\
&\mu_1(x_1)dx_1=\int_{x_0\in X}dM(x_0,x_1).\nonumber
\end{align}
The optimal mass transport cost is not necessarily a metric. However, 
if $X$ is a separable metric space with metric $d$ and we let $c(x_0,x_1)=d(x_0,x_1)^p$ where $p \ge 1$, then $T(\mu_0, \mu_1)^{1/p}$ is a metric on the set of nonnegative measures on $X$ with fixed mass.
This is the so called Wasserstein metrics.%
\footnote{The Wasserstein metric was, as many other concepts in this field, first defined by L. Kantorovich, and the naming is therefore somewhat controversial, see \cite{vershik2013long} and \cite[pp. 106-107]{villani2008optimal}.}
Moreover, $T(\mu_0, \mu_1)$ is weak$^*$ continuous on this set 
\cite[Theorem 6.9]{villani2008optimal}. Although the optimal mass transport is only defined for functions (measures) of the same mass, it can also be extended to handle measures with unbalanced masses \cite{georgiou2009metrics} (see also \cite{chizat2015unbalanced} for recent developments).

In this paper we consider the discrete version of the Kantorovich formulation \eqref{eq:Kantorovich}
\begin{align}
T(\mu_0, \mu_1) =\quad \min_{M\ge 0}\qquad  &\tr(C^T M)\nonumber\\
 \mbox{subject to}\quad &\mu_0=M \ett_{n_1}\label{eq:T}\\
& \mu_1=M^T \ett_{n_0}.\nonumber
\end{align}
In this setting the mass distributions are represented by two vectors $\mu_0\in \mR_+^{n_0}$ and $\mu_1\in \mR_+^{n_1}$, where the element $[\mu_k]_i$ corresponds to the mass in the point $x_{(k,i)}\in X$ for $i=1,\ldots, n_k$ and $k=0,1$. 
 A transference plan is represented by a matrix $M \in \mR_+^{n_0 \times n_1}$ where the value $m_{ij}:=[M]_{ij}$ denotes the amount of mass transported from point $x_{(0,i)}$ to $x_{(1,j)}$.
Such a plan is a feasible transference plan from 
$\mu_0$ to $\mu_1$ if the row sums of $M$ is $\mu_0$ and the column sums of $M$ is $\mu_1$.
The associated cost of a transference plan is $\sum_{i=1}^{n_0}\sum_{j=1}^{n_1}c_{ij}m_{ij}=\tr(C^TM)$, where $[C]_{ij}=c_{ij}=c(x_{(0,i)},x_{(1,j)})$ is the transportation cost from $x_{(0,i)}$ to $x_{(1,j)}$. 

Even though \eqref{eq:T} is a convex optimization problem, it is in many cases computationally infeasible due to the vast number of variables. The number of variables is $n_0n_1$ and if one seek to solve the optimal transport problem between two $256\times 256$ images, this results in  more than $4\cdot 10^{9}$ variables.

One recent method to approximate the optimal transport problem, and thereby allowing for addressing large problems, is to introduce an entropic regularizing term. This was proposed in \cite{cuturi2013sinkhorn} where the following optimization problem is considered
\begin{subequations} \label{eq:Teps}
\begin{align}
T_{\epsilon}(\mu_0,\mu_1)=\quad \min_{M\ge 0}\quad & \tr(C^T M)+ \epsilon D(M)\label{eq:Tepsa}\\
\mbox{ subject to }  & \mu_0=M \ett_{n_1}\label{eq:Tepsb}\\
&\mu_1=M^T \ett_{n_0},\label{eq:Tepsc}
\end{align}
\end{subequations}
where $D(M)=\sum_{i=1}^{n_0}\sum_{j=1}^{n_1} (m_{ij}\log(m_{ij})-m_{ij}+1)$ is a normalized entropy  term \cite{csiszar1991why}. This type of regularization is sometimes denoted entropic proximal  \cite{teboulle1992entropic}, and has previously been considered explicitly for linear programming \cite{fang1992unconstrained}. Also worth noting is that $D(M)$ is nonnegative and equal to zero if and only if $M = \ett_{n_0} \ett_{n_1}^T$.

A particularly nice feature with this problem is that any optimal solution belongs to an a priori known structure parameterized by $n_0+n_1$ variables via diagonal scaling (see \eqref{eq:structure} below). This can be seen  
by relaxing the equality constraints, and consider the Lagrange function
\begin{eqnarray}\label{eq:Teps_L}
L(M, \lambda_0, \lambda_1)&=&\tr(C^T M)+\epsilon D(M)\\&&+\lambda_0^T(\mu_0-M \ett_{n_1})+\lambda_1^T(\mu_1-M^T \ett_{n_0}).\nonumber
\end{eqnarray} 
For given dual variables, $\lambda_0\in \mR^{n_0}$ and $\lambda_1\in \mR^{n_1}$, the minimum $m_{ij}$ is obtained at
\begin{equation}\label{eq:Teps_grad}
0=\frac{\partial L(M,\lambda_0, \lambda_1)}{\partial m_{ij}}=c_{ij}+\epsilon \log(m_{ij})-\lambda_0(i)-\lambda_1(j)
\end{equation}
which can be expressed explicitly as $m_{ij}=e^{ \lambda_0(i)/\epsilon}e^{- c_{ij}/\epsilon }e^{\lambda_1(j)/\epsilon }$, or equivalently the solution is on the form
\begin{equation}
M=\diag(u_0)K \diag(u_1) \label{eq:structure}
\end{equation}
where $K=\exp(-C/\epsilon )$, $u_0=\exp(\lambda_0/\epsilon)$, and $u_1=\exp(\lambda_1/\epsilon)$. A theorem by Sinkhorn \cite{sinkhorn1967diagonal} states that for any matrix $K$ with positive elements, there exists diagonal matrices $\diag (u_0)$ and $\diag(u_1)$ with $u_0,u_1>0$ such that $M=\diag(u_0)K \diag(u_1)$ has prescribed row sums  and columns sums (i.e., $M$ satisfies \eqref{eq:Tepsb} and \eqref{eq:Tepsc}). Furthermore, the vectors $u_0$ and $u_1$ may be found by Sinkhorn iterations, i.e., alternatively solving \eqref{eq:Tepsb} for $u_0$ and \eqref{eq:Tepsc} for $u_1$: 
\begin{subequations} \label{eq:sinkhorn}
\begin{align}
&&\diag(u_0)K\diag(u_1) \ett =&\mu_0\quad\Rightarrow\quad u_0=\mu_0./(Ku_1)\label{eq:sinkhorn_a}\\
&&\diag(u_1)K^T\diag(u_0) \ett =&\mu_1\quad\Rightarrow\quad u_1=\mu_1./(K^Tu_0).\label{eq:sinkhorn_b}
\end{align}
\end{subequations}
The main computational bottlenecks in each iteration are the multiplications $Ku_1$ and $K^Tu_0$, and the iterations are therefore highly computationally efficient. In particular for cases where the matrix $K$ has a structure which can be exploited (see discussion in Section \ref{sec:prox}).
Furthermore the convergence rate is linear  \cite{franklin1989scaling} (cf. \cite{chen2016entropic} for generalization to positive functions).

Recently, in \cite{benamou2015iterative}, it was shown that the same iterative procedure for solving \eqref{eq:Teps} can also be recovered using Bregman projections \cite{bregman1967relaxation}, or Bregman-Dykstra iterations \cite{bauschke2000dykstras,boyle1986method}.
This approach was also used to solve several related problems such as for computing barycenters, multimarginal optimal transport problems and tomographic reconstruction \cite{benamou2015iterative}.
It has also been shown that the Sinkhorn iterations can be interpreted as block-coordinate ascent maximization of the dual problem via a generalization of the Bregman-Dykstra iterations \cite[Proposition~3.1]{peyre2015entropic}.
In Section~\ref{sec:dual} we derive this result using Lagrange duality. 
This observation opens up for enlarging the set of optimization problems that fit into this framework, and may also result in new algorithms adopted to the dual optimization problem.

\subsection{Variable splitting in convex optimization}\label{sec:splitting0}
Variable splitting is a technique for solving variational problems where the objective function is the sum of several terms that are simple in some sense (see, e.g., \cite{bauschke2011convex, combettes2011proximal, eckstein1989splitting}).
One of the more common algorithms for variable splitting is ADMM \cite{boyd2011distributed}, but there is a plentiful of other algorithms 
such as primal-dual forward-backward splitting algorithms \cite{bot2015convergence, chambolle2011first, sidky2012convex}, primal-dual forward-backward-forward splitting algorithms \cite{combettes2012primal}, and Douglas-Rachford type algorithms \cite{bot2013douglas, eckstein1992douglas}.
 For a good overview see \cite{komodakis2015playing}.
In this work we will explicitly consider variable splitting using a Douglas-Rachford type  algorithm presented in \cite{bot2013douglas}, but in order to better understand how the algorithm works we will first have a look at the proximal point algorithm and basic Douglas-Rachford variable splitting.

The basic idea behind many of the splitting techniques mentioned above spring from the so called proximal point algorithm for maximally monotone operators \cite{rockafellar1976monotone}. An operator $S : H \to H$, where $H$ is a real Hilbert space with the inner product $\langle \cdot, \cdot \rangle$, is called monotone if
\[
\langle z - z', w - w' \rangle \geq 0, \quad \text{for all } z,z' \in H, w \in S(z) \text{ and } w' \in S(z'), 
\]
and maximally monotone if in addition the graph of $S$,
\[
\{ (z, w) \in H \times H \mid w \in S(z) \},
\]
is not properly contained in the graph of any other monotone operator. The interest in such operators steam from the fact that the subdifferential of a proper, convex and lower semi-continuous function $f$, denoted by $\partial f$, is a maximally monotone operator \cite{rockafellar1976monotone}. Moreover, a global minimizer of such $f$ is any point $z$ so that $0 \in \partial f(z)$, which we denote $z \in \text{zer}(\partial f)$. 
For a maximally monotone operator $S$ and any scalar $\sigma > 0$ the operator  $(I +\sigma S)^{-1}$ is called the resolvent operator or proximal mapping. 
The proximal point algorithm is a fixed-point iteration of the resolvent operator,
\[
z^{k+1} = (I + \sigma S)^{-1} (z^k),
\]
and if $\text{zer}(S) \neq \emptyset$ then $z^k$ converges weakly to a point $z^\infty \in \text{zer}(S)$ \cite{rockafellar1976monotone}.  For the case that $S = \partial f$ the resolvent operator is called the proximal operator, and is given by 
\begin{equation}\label{eq:proximal}
(I + \sigma \partial f)^{-1}(z) = \Prox^{\sigma}_{f} (z) := \argmin_{z'} \left\{ f(z') + \frac{1}{2\sigma} \| z' - z \|_2^2 \right\}.
\end{equation}
Hence, fixed point iterations of the form
\[
z^{k+1} = \argmin_{z'} \left\{ f(z') + \frac{1}{2\sigma} \| z' - z^{k} \|_2^2 \right\}
\]
generates a sequence that converges weakly to a global minimizer of $f$. The parameter $\sigma$ determines the weighting between $f$ and the squared norm in $H$, and can be interpreted as a step length.

When the function to be minimized is a sum of several terms, then the resolvent operator of $S=A+B$, i.e., $(I + \sigma (A+B))^{-1}$, can be approximated in terms of the operators $A$ and $B$ and their resolvent operators $(I + \sigma A)^{-1}$ and $(I + \sigma B)^{-1}$ \cite{eckstein1989splitting}. 
This can give rise to fast schemes when the proximal operator of each term in the sum can be computed efficiently.  
One specific such algorithm, that is globally convergent, is the Douglas-Rachford splitting algorithm \cite{eckstein1992douglas}. In Section \ref{sec:inverse} we will use the splitting algorithm presented in \cite{bot2013douglas}, which extends this framework, in order to address a fairly general class of inverse problems.

\section{The dual problem and generalized Sinkhorn iterations}\label{sec:dual}

In this section we will see that the Sinkhorn iteration is identical to block-coordinate ascent of the  corresponding dual problem. Further, we will show that this procedure can also be applied to a set of inverse problems where the transportation cost is used as a regularizing term, and in particular for computing the proximal operator of $T_\epsilon(\mu_0,\cdot)$. 

\subsection{Sinkhorn iterations and the dual problem}

The Lagrangian dual of \eqref{eq:Teps} is defined as the minimum of \eqref{eq:Teps_L} over all $M\ge 0$ \cite[p. 485]{bertsekas1999nonlinear}. In our case this can 
be obtained by noting that  
\begin{eqnarray*}
&&\min_{M\ge 0} L(M, \lambda_0, \lambda_1)=L(M^*, \lambda_0, \lambda_1)\\
&&=\lambda_0^T\mu_0+\lambda_1^T\mu_1+\sum_{i=1}^{n_0}\sum_{j=1}^{n_1}\bigg(m_{ij}^*(c_{ij}
+\epsilon \log m_{ij}^*-\lambda_0(i)-\lambda_1(j)
)+1-m_{ij}^*\bigg)\\
&&=\lambda_0^T\mu_0+\lambda_1^T\mu_1+\epsilon\sum_{i=1}^{n_0}\sum_{j=1}^{n_1}(1-m_{ij}^*)\\
&&=\lambda_0^T\mu_0+\lambda_1^T\mu_1-\epsilon \exp(\lambda_0^T/\epsilon)\exp(-C/\epsilon)\exp(\lambda_1/\epsilon)+\epsilon n_0 n_1,
\end{eqnarray*}
where the optimal solution $M^*=[m^*_{ij}]_{ij}$  is specified by the equation \eqref{eq:Teps_grad}, which is also used in the third and forth equalities. This gives the following expression for the dual problem, a result which can also be found in \cite[Section 5]{cuturi2014fast}.
\begin{proposition}[{\cite{cuturi2014fast}}]\label{prp:dual_eps}
A Lagrange dual of \eqref{eq:Teps} is given by
\begin{align}
\max_{\lambda_0,\lambda_1} & \;\lambda_0^T\mu_0 +\lambda_1^T\mu_1-\epsilon \exp(\lambda_0^T/\epsilon)\exp(-C/\epsilon)\exp(\lambda_1/\epsilon)+\epsilon n_0 n_1.\label{eq:dual_eps}
\end{align}
\end{proposition}

Note the resemblance between the entropy relaxed dual formulation \eqref{eq:dual_eps} and the dual of the optimal transport problem \eqref{eq:T} \cite{villani2008optimal} 
\begin{subequations}\label{eq:T_dual}
\begin{align}
\max_{\lambda_0,\lambda_1} & \quad \mu_0^T\lambda_0+\mu_1^T\lambda_1\label{eq:T_duala}\\
\mbox{subject to } & \quad \lambda_0 \ett_{n_1}^T+\ett_{n_0} \lambda_1^T\le C.\label{eq:T_dualb}
\end{align}
\end{subequations}
The difference is that  the inequality constraint \eqref{eq:T_dualb} 
is exchanged for the penalty term
\begin{equation}\label{eq:barrier}
-\epsilon \exp(\lambda_0^T/\epsilon)\exp(-C/\epsilon)\exp(\lambda_1/\epsilon)
\end{equation}
in the objective function of \eqref{eq:dual_eps}. As $\epsilon\to 0$, the value of the barrier term \eqref{eq:barrier} goes to $0$ if the constraint \eqref{eq:T_dualb} is satisfied and to $-\infty$ otherwise.

Next, consider maximizing the dual objective \eqref{eq:dual_eps} with respect to $\lambda_0$ for a fixed $\lambda_1$.
This is attained by setting the gradient of the objective function in \eqref{eq:dual_eps} with respect to $\lambda_0$ equal to zero, hence $\lambda_0$ satisfies
\[
\mu_0= \exp(\lambda_0/\epsilon)\odot \left(\exp(-C/\epsilon)\exp(\lambda_1/\epsilon)\right).
\]
This is identical to the update formula \eqref{eq:sinkhorn_a} for $u_0=\exp(\lambda_0/\epsilon)$,  corresponding to the Sinkhorn iterations, where as before $u_1=\exp(\lambda_1/\epsilon)$. By symmetry, maximizing $\lambda_1$ for a fixed $\lambda_0$ gives a corresponding expression which is identical to \eqref{eq:sinkhorn_b}. Hence the Sinkhorn iterations corresponds to block-coordinate ascent in the dual problem, i.e., iteratively maximizing the objective in \eqref{eq:dual_eps} with respect to $\lambda_0$ while keeping $\lambda_{1}$ fixed, and vice versa.
\begin{corollary}[{\cite{peyre2015entropic}}] \label{cor:Sinkhorn}
The Sinkhorn iteration scheme \eqref{eq:sinkhorn} is a block-coordinate ascent algorithm of the dual problem \eqref{eq:dual_eps}. 
\end{corollary}

This was previously observed in \cite[Section 3.2]{peyre2015entropic}.
As we will see next, block-coordinate ascent of the dual problem results in fast Sinkhorn-type iterations for several different problems.

\subsection{Generalized Sinkhorn iterations}

Let us go back to the optimization problem \eqref{eq:Pgen0} that contains an optimal mass transport cost
\begin{align}\label{eq:Pgen1}
\min_{\mue} &\quad  T_\epsilon(\mu_0, \mue)+g(\mue), 
\end{align}
where $\mu_0$ is a prior, and $g$ is a term that could include other regularization terms and data miss-match. In order to guarantee that this problem has a solution and is convex, we introduce the following assumption.  

\begin{assumption}\label{ass:g}
Let $g$ be a proper, convex and lower semi-continuous function that is finite in at least one point with mass equal to $\mu_0$, i.e., $g(\mue)<\infty$ for some $\mue$ with $\sum_{i=1}^{n_0} \mu_0(i) = \sum_{j=1}^{n_1} \mue(j)$.
\end{assumption}

The first part of this assumption is to make the problem convex, and the second part is imposed so that \eqref{eq:Pgen1} has a feasible solution. Moreover, note that $T_\epsilon(\mu_0, \mue)<\infty$ restricts $\mue$ to a compact set, which guarantees the existence of an optimal solution.
Using the definition of the optimal transport cost, the problem \eqref{eq:Pgen1} can equivalently be formulated as
\begin{align}
\min_{M\ge 0,  \mue} & \tr(C^TM)+\epsilon D(M)+g(\mue)\label{eq:Pgen_eps0}\\
\mbox{subject to } &  \mu_0=M \ett_{n_1}\nonumber\\
& \mue=M^T \ett_{n_0}.\nonumber
\end{align}
The Lagrangian dual problem of \eqref{eq:Pgen_eps0} can be obtained using the same steps as the derivation of Proposition~\ref{prp:dual_eps}. See Appendix~\ref{app:A} for details. Results similar to Proposition \ref{prp:dual} is also obtained in the recent preprints \cite[Theorem 1]{chizat2016scaling} and \cite{schmitzer2016stabilized}.
 
\begin{proposition}\label{prp:dual}
Let $\mu_0>0$ be given and let $g$ satisfy Assumption~\ref{ass:g}. Then the Lagrange dual of \eqref{eq:Pgen_eps0} is given by
\begin{align}
\max_{\lambda_0,\lambda_1} & \;\lambda_0^T\mu_0 -g^*(-\lambda_1)-\epsilon \exp(\lambda_0^T/\epsilon)\exp(-C/\epsilon)\exp(\lambda_1/\epsilon)+\epsilon n_0 n_1\label{eq:dualgen_eps}
\end{align}
and strong duality holds.
\end{proposition}

The only difference between \eqref{eq:dualgen_eps} and \eqref{eq:dual_eps} is that the term $\lambda_1^T\mu_1$ is replaced by $-g^*(-\lambda_1)$, where $g^*$ denotes the dual (or Fenchel) conjugate functional
\[
g^*(\lambda)=\sup_{\mu} \left(\lambda^T \mu-g(\mu)\right).
\]
Clearly, Proposition~\ref{prp:dual_eps} is a special case of Proposition~\ref{prp:dual}, and the optimization problem $T_{\epsilon}(\mu_0,\mu_1)$ in \eqref{eq:Teps} is recovered from \eqref{eq:Pgen1} 
if $\mue$ is fixed to $\mu_1$, i.e., $g(\mue)=\cI_{\mu_1}(\mue)$. 
Since \eqref{eq:dualgen_eps} is a dual problem, the objective function is concave \cite[p. 486]{bertsekas1999nonlinear}, but not necessarily strictly concave (e.g., as in the case \eqref{eq:dual_eps}).
Moreover, Assumption~\ref{ass:g} assures strong duality between \eqref{eq:Pgen_eps0} and \eqref{eq:dualgen_eps} (see proof in Appendix~\ref{app:A}).
As for the standard optimal mass transport problem, we now consider block-coordinate ascent to compute an optimal solution. The corresponding optimality conditions are given in the following lemma and are obtained by noting that the optimum is only achieved when zero is a (sub)gradient of \eqref{eq:dualgen_eps} \cite[pp. 711-712]{bertsekas1999nonlinear}.    
\begin{lemma}\label{lm:grad_cond}
For a fixed $\lambda_1$, then $\lambda_0$ is the maximizing vector of \eqref{eq:dualgen_eps} if 
\begin{subequations}\label{eq:dual_stat2}
\begin{equation}
\mu_0= \exp(\lambda_0/\epsilon)\odot \left(\exp(-C/\epsilon)\exp(\lambda_1/\epsilon)\right).\label{eq:dual_stat2a}
\end{equation}
Similarly, for a fixed $\lambda_0$, then $\lambda_1$ is the maximizing vector of \eqref{eq:dualgen_eps} if 
\begin{equation}\label{eq:dual_stat2b}
0\in \partial g^*(-\lambda_1) -  \exp(\lambda_1/\epsilon)  \odot  \left( \exp(-C^T/\epsilon)\exp(\lambda_0/\epsilon)\right).
\end{equation}
\end{subequations}
\end{lemma}

Whenever \eqref{eq:dual_stat2} can be computed efficiently, block-coordinate ascent could give rise to a fast computational method for solving \eqref{eq:dualgen_eps}. In particular this is the case if the equation \eqref{eq:dual_stat2b} can be solved element by element, i.e., in $\mathcal{O}(n)$ (excluding the matrix-vector multiplication $\exp(-C^T/\epsilon)\exp(\lambda_0/\epsilon)$). This is true for several cases of interest. 

\begin{example}
$g(\mu)=\cI_{\{\mu_1\}}(\mu)$. This corresponds to the optimal mass transport problem \eqref{eq:Teps} 
\end{example}

\begin{example}
$g(\mu)=\|\mu-\mu_1\|_1$.
\end{example}

\begin{example}\label{ex:prox}
$g(\mu)=\frac{1}{2}\|\mu-\mu_1\|_2^2$.
\end{example}

\begin{example}
$g(\mu)=\frac{1}{2}\|A\mu-\mu_1\|_2^2$ where  $A^*A$ is diagonal and invertible.
\end{example}

The example~\ref{ex:prox} is of particular importance, since this corresponds to computing the proximal operator of the transportation cost, and will be addressed in detail in the next subsection. 
Note that Lemma~\ref{lm:grad_cond} can also be expressed in terms of the entropic proximal operator \cite{peyre2015entropic}, see Appendix~\ref{app:C} for details.

\subsection{Sinkhorn-type iterations for evaluating the proximal operator}\label{sec:prox}
The proximal point algorithm \cite{rockafellar1976monotone} and splitting methods \cite{eckstein1989splitting} are extensively used in optimization, and a key tool is the computation of the proximal operator \eqref{eq:proximal} (see Section~\ref{sec:splitting0} and Section~\ref{sec:inverse}).
In order to use this kind of methods for solving problems of the form \eqref{eq:Pgen1}
we consider the proximal operator of the entropy regularized mass transport cost, and propose a Sinkhorn-type algorithm for computing the proximal operator of the transportation cost
\begin{equation*}
\Prox^\sigma_{ T_\epsilon(\mu_0, \cdot)}(\mu_1)=\argmin_{\mue} T_\epsilon(\mu_0,\mue)+ \frac{1}{2\sigma}\|\mue-\mu_1\|_2^2.
\end{equation*} 
First note that this can be identified with the optimization problem \eqref{eq:Pgen_eps0} where the data fitting term and the corresponding conjugate functional are 
\begin{align*}
g(\mu)&=\frac{1}{2\sigma}\|\mu-\mu_1\|^2_2,\quad & g^*(\lambda)&=\lambda^T\left(\mu_1+\frac{\sigma}{2}\lambda\right).
\end{align*}
The dual problem is then given by
\begin{align}\label{eq:dual_prox}
\max_{\lambda_0,\lambda_1} & \;\lambda_0^T\mu_0 +\lambda_1^T(\mu_1-\frac{\sigma}{2}\lambda_1)-\epsilon \exp(\lambda_0^T/\epsilon)\exp(-C/\epsilon)\exp(\lambda_1/\epsilon)+\epsilon n_0 n_1.
\end{align}
The optimality conditions corresponding to Lemma~\ref{lm:grad_cond} leads to the equations (see the proof of Theorem \ref{thm:main} in Appendix~\ref{app:B} for the derivation)
\begin{subequations}\label{eq:equations_prox}
\begin{align}
\lambda_0 &= \epsilon \log\left( \mu_0./(K\exp(\lambda_1/\epsilon))\right)\label{eq:dual_sink}\\
\lambda_1 &=\frac{\mu_1}{\sigma}-\epsilon\omega\left(\frac{\mu_1}{\sigma\epsilon}+\log\left(K^T\exp(\lambda_0/\epsilon)\right)-\log(\sigma\epsilon)\right).\label{eq:dual_prox2}
\end{align}
\end{subequations}
Here $\omega$ denotes the elementwise Wright $\omega$ function,%
\footnote{Our implementation use $\omega(x)=W(e^x)$ for $x\in \mR$, where $W$ is the Lambert $W$ function \cite{corless2002wright}.}
 i.e., the function mapping $x\in\mR$ to $\omega(x)\in \mR_+$ for which $x=\log(\omega(x))+\omega(x)$ \cite{corless2002wright}.
The first equation can be identified with the first update equation \eqref{eq:sinkhorn_a} in the Sinkhorn iteration.  Note that the bottlenecks in the iterations \eqref{eq:equations_prox} are the multiplications with of $K$ and $K^T$. All other operations are elementwise and can hence be computed in ${\mathcal O}(n)$ where $n=\max(n_0,n_1)$. The full algorithm is presented in Algorithm~\ref{alg:prox_Sinkhorn}. This leads to one of our main results.

\begin{theorem}\label{thm:main}
The variables $(\lambda_0, \lambda_1)$ in Algorithm~\ref{alg:prox_Sinkhorn} converges to the optimal solution of the dual problem \eqref{eq:dual_prox}. Furthermore, the convergence rate is locally q-linear.
\end{theorem}

\noindent{\bf Proof:} See Appendix~\ref{app:B}.\hfill$\square$\\[-10pt] 

This proof is based on the duality (i.e., Proposition \ref{prp:dual} and Lemma \ref{lm:grad_cond}). The algorithm could also be derived directly using Bregman projections \cite{bregman1967relaxation}, similarly to the derivation of the Sinkhorn iteration in \cite{benamou2015iterative}. Some remarks are in order regarding the computation of the iterations for the Sinkhorn-type iterations. 

\begin{remark}\label{rem:Fast_K}
The bottlenecks in the iterations \eqref{eq:equations_prox} are the multiplications with the matrices $K$ and $K^T$. All other operations are elementwise. In many cases of interest the structures of $K$ can be exploited for fast computations. In particular when the mass points (pixel/voxel locations) $x_{(0,i)}=x_{(1,i)}$ are on a regular grid and the cost function is translation invariant, e.g., as in our application  example (see \eqref{eq:costfucntion}) where the cost function only depends on the distance between the grid points.
Then the matrices $C$ and $K$ are multilevel Toeplitz-block-Toeplitz and the multiplication can be performed in ${\mathcal O}(n\log(n))$ using the fast Fourier transform (FFT) (see, e.g., \cite{lee1986fast}).
\end{remark}

\begin{remark}
In order for \eqref{eq:Teps} to approximate the optimal mass transportation problem \eqref{eq:T} it is desirable to use a small $\epsilon$. 
The entropy regularization has a smoothing effect on the transference plan $M$ and a too large value of $\epsilon$ may thus result in an undesirable solution to the variational problem \eqref{eq:Pgen1}.
However, as $\epsilon\to 0$ the problem becomes increasingly ill-conditioned and the convergence becomes slower \cite{cuturi2013sinkhorn}. To handle the ill-conditioning one can 
stabilize the computations using logarithmic reparameterizations.
One such approach is described in \cite{chen2016entropic}, where the variables $\log(u_0), \log(u_1)$ are used, together with appropriate normalization and truncation of the variables, to compute the Sinkhorn iterations \eqref{eq:sinkhorn}.
Another approach to handle the ill-conditioning is described in \cite{schmitzer2016stabilized}, where a different logarithmic reparameterization is used together with an adaptive scheme for scaling of both $\epsilon$ and the discretization grid (cf. \cite{chizat2016scaling}).
However, note that these approaches are not compatible with utilizing FFT-computations for the matrix-vector products by exploiting the Toeplitz-block-Toeplitz structure in $C$ (and $K$). Due to this fact we have not used this type of stabilization.
\end{remark}

\renewcommand{\algorithmicrequire}{\textbf{Input:}}
\renewcommand{\algorithmicensure}{\textbf{Output:}}

\algsetup{indent=12pt}

\begin{algorithm}
\caption{Generalized Sinkhorn algorithm for evaluating the proximal operator of $T_\epsilon(\mu_0,\cdot)$.}
\label{alg:prox_Sinkhorn}
\begin{algorithmic}[1]
\REQUIRE $\epsilon$, $C$, $\lambda_0$, $\mu_0$, $\mu_1$
\STATE $K=\exp(-C/\epsilon)$
\WHILE{Not converged}
\STATE $\lambda_0 \leftarrow \epsilon \log\left( \mu_0./(K\exp(\lambda_1/\epsilon))\right)$\label{alg:step_sinkhorn}
\STATE $\lambda_1 \leftarrow \frac{\mu_1}{\sigma}-\epsilon\omega\left(\frac{\mu_1}{\sigma\epsilon}+\log\left(K^T\exp(\lambda_0/\epsilon)\right)-\log(\sigma\epsilon)\right)$\label{alg:step_prox}
\ENDWHILE
\ENSURE $\mue \leftarrow  \exp(\lambda_1/\epsilon) \odot (K^T \exp(\lambda_0/\epsilon))$
\end{algorithmic}
\end{algorithm}

\section{Inverse problems with optimal mass transport priors}\label{sec:inverse}
In this section we use the splitting framework \cite{bot2013douglas} to formulate and solve inverse problems. This is a generalization of the Douglas-Rachford algorithm \cite{eckstein1992douglas}. In particular this allows us to address large scale problems of the form \eqref{eq:Pgen1}, since the proximal operator of the regularized transport problem can be computed efficiently using generalized Sinkhorn iterations.

The theory in \cite{bot2013douglas} provides a general framework that allows for solving a large set of convex optimization problems. Here we consider problems of the form 
\begin{equation}\label{eq:general_opt_prob}
\inf_{z \in H} f(z) + \sum_{i = 1}^m g_i(L_iz - r_i)
\end{equation}
where $H$ is a Hilbert space, $f : H \rightarrow \bar{\mathbb{R}}$ is proper, convex and lower semi-continuous, $g_i : G_i \rightarrow \bar{\mathbb{R}}$, where $G_i$ is a Hilbert space and $g_i$ is is proper, convex and lower semi-continuous,  and $L_i : H \rightarrow G_i$ is a nonzero bounded linear operator, which is a special case of the structure considered in \cite{bot2013douglas}.

The problem \eqref{eq:general_opt_prob} can be solved by the iterative algorithm \cite[(3.6)]{bot2013douglas}, in which we only need to evaluate the proximal operators $\Prox^\tau_f$ and $\Prox^{\sigma_i}_{g_i^*}$ for $i = 1, \ldots, m$. For reference, this simplified version of the algorithm is shown in Algorithm~\ref{alg:DR}.
Note that by Moreau decomposition we have that $\Prox^\tau_f(x) = x - \tau \Prox^{1/\tau}_{f^*}(x/\tau)$ \cite[Theorem 14.3]{bauschke2011convex}, and therefore Algorithm~\ref{alg:DR} can be applied as long as either the proximal operators of the functionals $f$ and $\{g_i\}_{i=1}^m$ or the proximal operators of their Fenchel conjugates can be evaluated in an efficient way.

In the following examples we restrict the discussion to the  finite dimensional setting where the underlying set $X$ is a $d$-dimensional  regular rectangular grid with points $x_i$ for $i=1,\dots n$, hence the corresponding Hilbert space is $H=\mR^n$. Let $A:\mR^n\to \mR^m$ be a linear operator (matrix) representing measurements and let $\nabla:\mR^n\to \mR^{d\times n}$ be a discrete gradient operator%
\footnote{In each dimension of the regular rectangular grid a forward-difference is applied and the boundary of the domain is padded with zeros, i.e., $(\nabla \mu)_{j,i}$ is the forward difference along the $j$-axis in the grid point $x_i$.  See Appendix~\ref{app:D} for more details.}
based on the grid $X$. 
Furthermore, for $Y=(y_1,\ldots, y_n)\in \mR^{d\times n}$ we let $\|Y\|_{2,1}=\sum_{i=1}^n \|y_i\|_2$ be the isotropic $\ell_1$-norm (sometimes called group $\ell_1$-norm). With this notation $\|\nabla\mu\|_{2,1}$ is the isotropic total variation (TV) of $\mu$, which is often used as a convex surrogate for the support of the gradient \cite{rudin1992nonlinear, candes2005signal}. 
This terminology allows us to set up a series of optimization problem for addressing inverse problems.

\begin{algorithm}
\caption{Douglas-Rachford type primal-dual algorithm \cite{bot2013douglas}.}
\label{alg:DR}
\begin{algorithmic}[1]
\REQUIRE $\tau, (\sigma_i)_{i = 1}^m, (\lambda_n)_{n \geq 1}$, such that $\sum_{n = 1}^\infty \lambda_n(2 - \lambda_n)\! =\! \infty$ and $\tau \sum_{i = 1}^m \sigma_i \| L_i \|^2\! < 4$.
\STATE $n = 0$
\STATE $p_{1,0} = w_{1,0} = z_{1,0} = 0$
\STATE $p_{2,i,0} = w_{2,i,0} = z_{2,i,0} = 0$
\WHILE{Not converged}
\STATE $n \leftarrow n+1$
\STATE $p_{1,n} \leftarrow \Prox^\tau_f \left( x_n - \frac{\tau}{2}\sum_{i = 1}^m L_i^*v_{i,n} \right)$
\STATE $w_{1,n} \leftarrow 2p_{1,n} - x_n$
\FOR{i = 1, \ldots, m}
\STATE $p_{2,i,n} \leftarrow \Prox^{\sigma_i}_{g^*_i} \left( v_{i,n} + \frac{\sigma_i}{2}L_i w_{1,n} - \sigma_i r_i\right)$
\STATE $w_{2,i,n} \leftarrow 2p_{2,i,n} - v_{i,n}$
\ENDFOR
\STATE $z_{1,n} \leftarrow w_{1,n} - \frac{\tau}{2}\sum_{i = 1}^m L_i^* w_{2,i,n}$
\STATE $x_{n+1} \leftarrow x_n + \lambda_n(z_{1,n} - p_{1,n})$
\FOR{i = 1, \ldots, m}
\STATE $z_{2,i,n} \leftarrow w_{2,i,n} + \frac{\sigma_i}{2}L_i(2z_{1,n}-w_{1,n})$
\STATE $v_{i, n+1} \leftarrow v_{i,n} + \lambda_n(2z_{2,i,n} - p_{2,i,n})$
\ENDFOR
\ENDWHILE
\ENSURE $x^* = \Prox^\tau_f \left( x_n - \frac{\tau}{2}\sum_{i = 1}^m L_i^* v_{i,n} \right)$
\end{algorithmic}
\end{algorithm}

\subsection{Optimal mass transport priors using variable splitting}
Next, we use optimal mass transport for incorporating a priori information. We consider a particular case of \eqref{eq:Pgen1} for achieving this, and formulate the reconstruction problem as follows:
\begin{align}\label{eq:OMT}
\min_\mue \qquad &\gamma T_\epsilon (\mu_0, \mue) + \| \nabla \mue \|_{2,1} \\
\mbox{subject to }\quad &\|A\mue - b\|_2\le \kappa.\nonumber
\end{align}
Here $\mu_0$ is a prior, $\kappa$ quantifies the allowed measurement error, and $\gamma$ determines the trade off between the optimal transport prior and the TV-regularization.  
Since we can compute the proximal operator of $T_\epsilon (\mu_0, \mu)$ in an efficient way, this problem can be solved by, e.g., making a variable splitting according to
\[
\begin{array}{lll}
f(\cdot) = \gamma T_\epsilon (\mu_0, \cdot), \\
g_1(\cdot) = \| \cdot \|_{2,1}, & L_1 = \nabla, & r_1=0, \\
g_2(\cdot) = \cI_{{\cal B}_m(\kappa)}(\cdot), & L_2 = A,  & r_2=b, \\
\end{array}
\]
where ${\cal B}_m(\kappa)=\{\hat b\in \mR^m\,:\,\|\hat b\|_2\le \kappa\}$ is the  $\kappa$-ball in $\mR^m$. We apply Algorithm~\ref{alg:DR} for solving this problem and use Algorithm~\ref{alg:prox_Sinkhorn} for computing the proximal operator of $T_\epsilon (\mu_0, \cdot)$. 
Explicit expression for the proximal operator of the Fenchel dual of $g_1$ and $g_2$ can be computed, see, e.g., \cite{komodakis2015playing} for details.

\begin{remark} \label{rem:otherTom2}
An alternative first-order method to solve optimization problems with one or more transportation costs is considered in \cite{cuturi2015asmoothed}. 
The authors use a dual forward-backward (proximal-gradient) scheme where a key component is the evaluation 
of the gradient of the dual conjugate functional of $T_\epsilon(\mu_0,\cdot)$. However, our problem \eqref{eq:OMT} contains two terms in addition to the optimal mass transport cost and  does not directly fit into this framework. Since our method builds on the splitting framework in \cite{bot2013douglas} it allows for an arbitrary number of cost terms (see \eqref{eq:general_opt_prob}).
\end{remark}

\subsection{Formulating standard inverse problems using variable splitting}
Given a linear forward operator $A$, many common regularization methods for solving inverse problems can be formulated as optimization problems on the form \eqref{eq:general_opt_prob}. Hence they can also be solved using variable splitting. We will use this to compare the proposed reconstruction method with two other approaches. First with an approach using TV-regularization \cite{candes2005signal}, and second with an approach where we use a prior information with respect to the standard $\ell_2$-norm. We formulate and solve both these problems using
Algorithm~\ref{alg:DR}.

First, we consider the TV-regularization problem 
\begin{align}\label{eq:TV}
\min_{\mue\ge 0} \qquad & \| \nabla \mue \|_{2,1} \\
\mbox{subject to }\quad &\|A\mue - b\|_2\le \kappa \nonumber
\end{align}
and formulate it in the setting of \eqref{eq:general_opt_prob} by defining
\[
\begin{array}{lll}
f(\cdot) = \cI_{\mR_+^n}(\cdot), \\
g_1(\cdot) = \| \cdot \|_{2,1}, & L_1 = \nabla, & r_1=0, \\
g_2(\cdot) = \cI_{{\cal B}_m(\kappa)}(\cdot), & L_2 = A, & r_2=b.
\end{array}
\]
The positivity constraint is handled by $f$, the TV-regularization by $g_1$ and the data matching by $g_2$. 
All functions needed in Algorithm~\ref{alg:DR} can be explicitly computed \cite{komodakis2015playing}.

In \eqref{eq:TV} there is no explicit notion of prior information and reconstruction is entirely based on the data and an implicit assumption on the sparsity of the gradient. 
One way to explicitly incorporate prior information in the problem is to add an $\ell_2$-norm term that penalizes deviations from the given prior. This leads to the optimization problem
\begin{align}\label{eq:L2}
\min_{\mue\ge 0} \qquad & \gamma\|\mue-\mu_0\|_2^2 +\| \nabla \mue \|_{2,1} \\
\mbox{subject to }\quad &\|A\mue - b\|_2\le \kappa\nonumber
\end{align}
which we formulate it in the setting of \eqref{eq:general_opt_prob} by defining
\[
\begin{array}{lll}
f(\cdot) = \cI_{\mR_+^n}(\cdot), \\
g_1(\cdot) = \gamma \| \cdot \|_2^2, & L_1 = I, & r_1=\mu_0  \\
g_2(\cdot) = \| \cdot \|_{2,1}, & L_2 = \nabla, & r_2=0 \\
g_3(\cdot) = \cI_{{\cal B}_m(\kappa)}(\cdot), & L_3 = A, & r_3=b.
\end{array}
\]
The positivity constraint is handled by $f$, the $\ell_2$ prior by $g_1$, the TV-regularization by $g_2$ and the data matching by $g_3$.
The parameter $\gamma$ determines the trade off between the $\ell_2$ prior and the TV-regularization. 
Also here, all functions needed in Algorithm~\ref{alg:DR} can be computed explicitly \cite{komodakis2015playing}.

\section{Application in Computerized Tomography}\label{sec:applications}

\emph{Computerized Tomography} (CT) is an imaging modality that is frequently used in many areas, especially in medical imaging (see, e.g.,  the monographs \cite{bertero1998introduction, kak2001principles, natterer2001themathematics, natterer2001mathematical}). In CT the object is probed with X-rays, and since different materials attenuate X-rays to different degrees, the intensities of the incoming and outgoing X-rays contain information of the material content and distribution. In the simplest case, where the attenuation is assumed to be energy independent and where scatter and nonlinear effects are ignored, one gets the equation \cite[Chapter~3]{natterer2001mathematical}
\begin{equation*}
\int_L \mut(x) dx = \log \left( \frac{I_0}{I} \right).
\end{equation*}
Here $\mut(x)$ is the attenuation in the point $x$, $L$ is the line along which the X-ray beam travels through the object, and $I_0$ and $I$ are the the incoming and outgoing intensities. By taking several measurements along different lines $L$, one seek to  reconstruct the attenuation map $\mut$.

A set of measurements thus corresponds to the line integral of $\mut$ along a limited set of lines,  
and the corresponding operator that map $\mut$ to the line integrals is called a \emph{ray transform} or a \emph{partial Radon transform}. Let $A$ be the partial Radon transform operator, i.e., the operator such that $A(\mu)$ gives the line integral of $\mu$ along certain lines. This is a linear operator, and we consider the inverse problem of recovering $\mut$ from measurements 
\[
b = A(\mut) + \text{noise}.
\]
However, this is an ill-posed  inverse problem \cite[p. 40]{engl2000regularization}. In particular, the problem is severely ill-posed if the set of measurements is small or limited to certain angles, 
and hence regularization is needed to obtain an estimate $\mue$ of $\mut$.
One way to obtain such a $\mue$ is to formulate variational problems akin to the ones in Section~\ref{sec:inverse}.
In this section we consider computerized tomography problems, such as image reconstruction from limited-angle measurements, and use optimal mass transport to incorporate prior information to compensate for missing measurements. We also compare this method with standard reconstruction techniques.
Tomography problems with transport priors have previously been considered in \cite{abraham2017tomographic, benamou2015iterative}, but in a less general setting. We compare and discuss the details in Remark~\ref{rem:otherTom} in the end of this section.

\begin{figure}%
  \centering
  \hfil
  \subfloat[\label{fig:phantom}]{\includegraphics[width=0.5\textwidth]{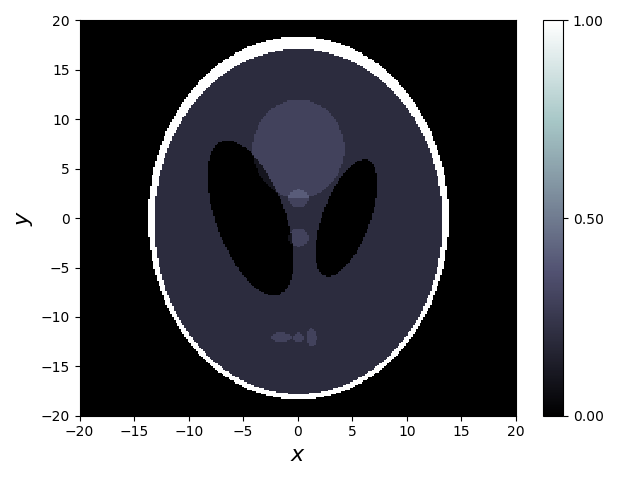}}%
  \hfil
  \subfloat[\label{fig:prior}]{\includegraphics[width=0.5\textwidth]{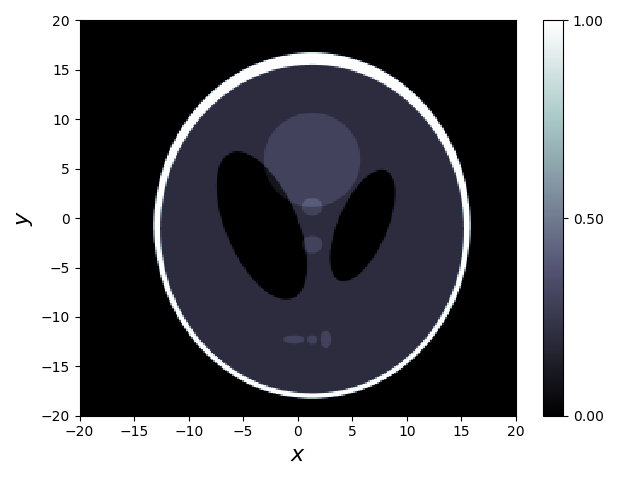}}%
  \caption{Figure showing (a) the Shepp-Logan phantom, (b) the deformed Shepp-Logan prior used in the first example.  
  Gray scale values are shown to the right of each image.}%
  \label{fig:phantom_and_prior}%
\end{figure}

\subsection{Numerical simulations}
\begin{figure}%
  \centering
  \hfil
  \subfloat[\label{fig:phantomHand}]{\includegraphics[width=0.5\textwidth]{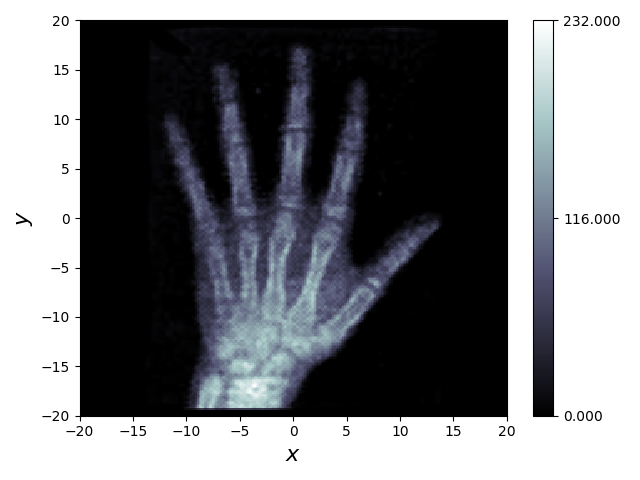}}%
  \hfil
  \subfloat[\label{fig:priorHand}]{\includegraphics[width=0.5\textwidth]{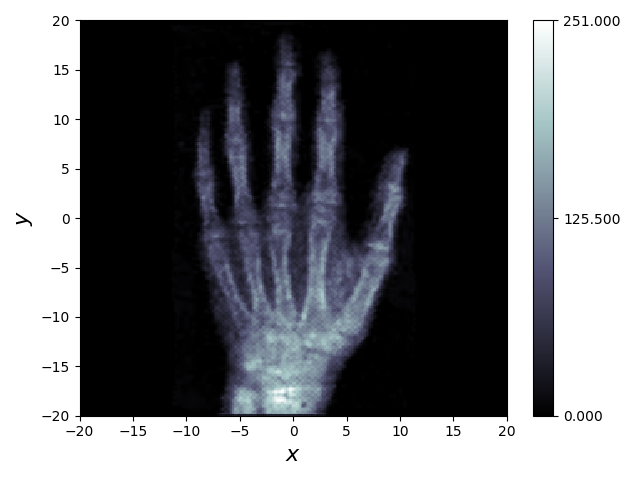}}%
  \caption{Figure showing (a) the hand image used as phantom, and (b) the hand image used as prior in the second example. Gray scale values are shown to the right of each image.}%
  \label{fig:phantom_and_priorHand}%
\end{figure}

 To this end, consider the Shepp-Logan phantom in Figure~\ref{fig:phantom} \cite{shepp1974fourier}, and the hand image in Figure~\ref{fig:phantomHand} \cite{bauer2015diffeomorphic}. 
 Assume that the deformed images in Figure~\ref{fig:prior} and Figure~\ref{fig:priorHand} has been reconstructed previously from a detailed CT scan of the patient (where the deformation is due to, e.g., motion or breathing). 
By using the deformed images \ref{fig:prior} and \ref{fig:priorHand} as prior information we want to reconstruct the images \ref{fig:phantom} and \ref{fig:phantomHand} from relatively few measurements of the later. For the reconstruction of the Shepp-Logan image we consider a scenario where the set of angles belong to a limited interval (see next paragraph). For the reconstruction of the hand image phantom we consider uniform spacing across all angles.   
In both examples we compare the reconstructions obtained by solving the variational problems \eqref{eq:OMT}, \eqref{eq:TV}, and \eqref{eq:L2}, as well as a standard filtered backprojection reconstruction \cite{natterer2001themathematics}.

\def \imDim {$256$}
\def \imDimTot {$65 \, 536$} 
\def \noiseLev {$5\%$}
\def \noiseLevHand {$3\%$}
\def \Angles {$30$}
\def \AnglesHand {$15$}
\def \AngleOne {$\pi/4$}
\def \AngleTwo {$3\pi/4$}
\def \AngleOneHand {$0$}
\def \AngleTwoHand {$\pi$}
\def \parLines {$350$}
\def \nrMeasurements {$10 \, 500$}
\def \nrMeasurementsHand {$5 \, 250$}
\def \nrVariables {$4 \cdot 10^9$} 
\def \DRiter {$10 \, 000$}
\def \SinkhornIter {$200$}

In these examples, the images  
have a resolution of \imDim $\times$\imDim\ pixels and we compute the data from the phantoms in Figures \ref{fig:phantom} and \ref{fig:phantomHand}. 
For the Shepp-Logan example we let data be collected from \Angles\ equidistant angles in the interval $[$\AngleOne, \AngleTwo $]$, and for the hand example from \AnglesHand\ equidistant angles in the interval $[$\AngleOneHand, \AngleTwoHand $]$. In both cases the data is the line integrals from \parLines\ parallel lines for each angle. 
On each data set white Gaussian noise is added (\noiseLev\, and \noiseLevHand, respectively). 
The corresponding optimization problems \eqref{eq:OMT}, \eqref{eq:TV}, and \eqref{eq:L2} are 
solved\footnote{A fixed number of \DRiter\, Douglas-Rachford iterations are computed in each reconstruction.} with the Douglas-Rachford type algorithm (Algorithm~\ref{alg:DR}) from \cite{bot2013douglas} where the functions are split according to the description in Section~\ref{sec:inverse}. 
The exact parameter values are given in Appendix~\ref{app:D}.
Since the grid is regular and the cost in the transportation cost is spatially invariant we can use FFT for the computations of $T_\epsilon(\mu_0,\cdot)$ and the corresponding proximal operator, and thus we do not need to explicitly store the cost matrix or the final transference plan (see Remark~\ref{rem:Fast_K} and Remark \ref{rem:big}).
These examples have been implemented and solved using ODL%
\footnote{Open source code, available at \url{https://github.com/odlgroup/odl}}  \cite{adler2017ODL},
which is a python library for fast prototyping focusing on inverse problems. The ray transform computations are performed by the GPU-accelerated version of ASTRA%
\footnote{Open source code, available at \url{https://github.com/astra-toolbox/astra-toolbox}}  \cite{palenstijn2011performance,vanaerle2015astra}.  
The code used for these examples is available online at \url{http://www.math.kth.se/~aringh/Research/research.html}.

\begin{figure}%
  \centering
  \hfil
  \subfloat[\label{fig:fbp}]{\includegraphics[width=0.5\textwidth]{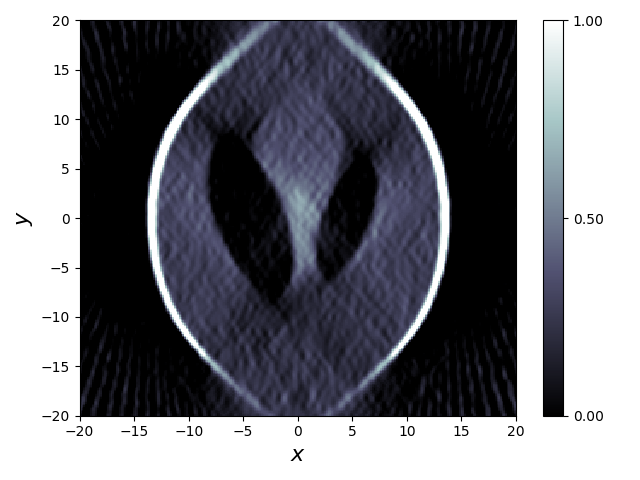}}%
  \hfil
  \subfloat[\label{fig:TV}]{\includegraphics[width=0.5\textwidth]{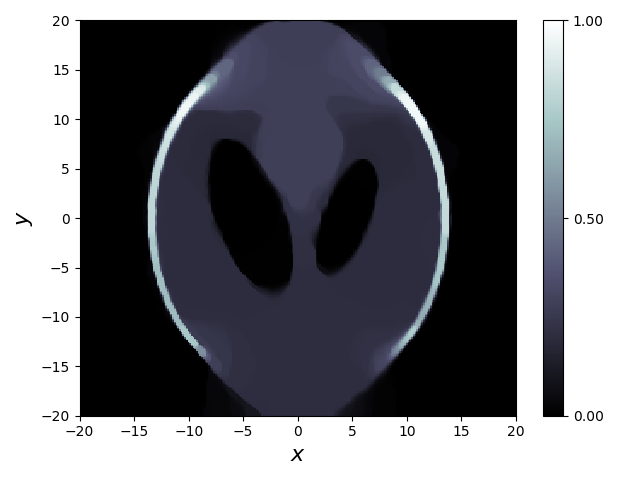}}%
  \hfil
  \subfloat[\label{fig:L2}]{\includegraphics[width=0.5\textwidth]{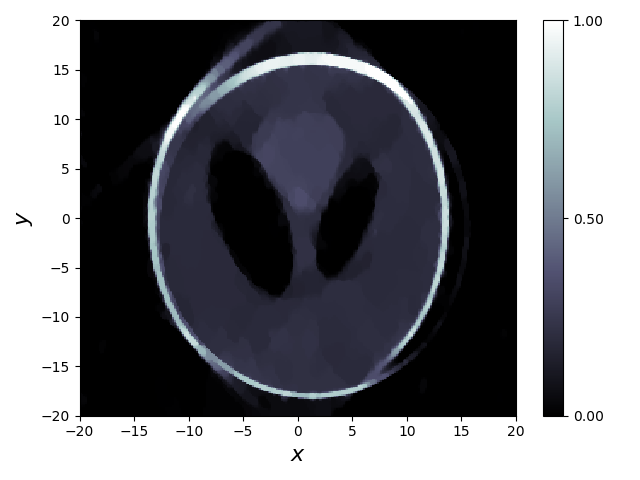}}%
  \hfil
  \subfloat[\label{fig:OMT}]{\includegraphics[width=0.5\textwidth]{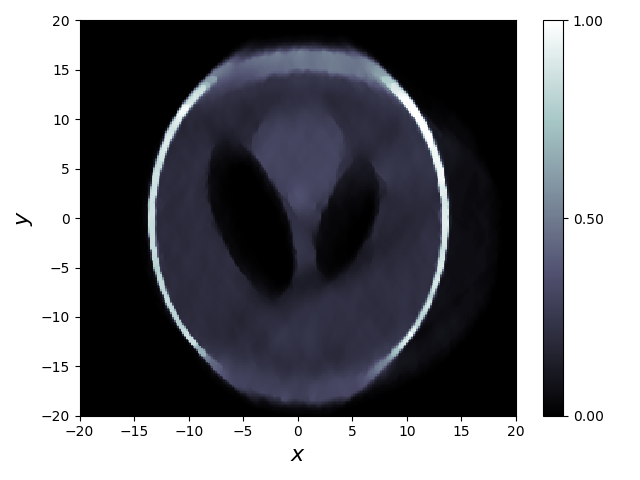}}%
  \caption{Reconstructions using different methods. (a) Filtered backprojection, 
 (b) reconstruction using TV-regularization, (c) reconstruction with $\ell_2^2$-prior and TV-regularization ($\gamma=10$), and (d) reconstruction with optimal transport prior and TV-regularization ($\gamma=4$).}%
  \label{fig:results}%
\end{figure}

\begin{figure}%
  \centering
  \hfil
  \subfloat[\label{fig:L2_a}]{\includegraphics[width=0.33\textwidth]{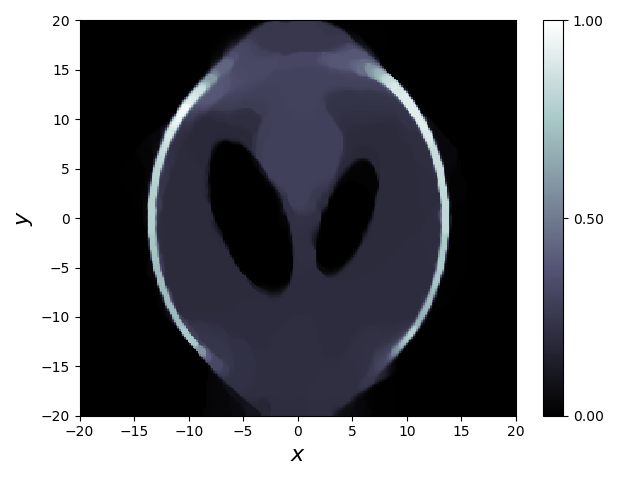}}%
  \subfloat[\label{fig:L2_b}]{\includegraphics[width=0.33\textwidth]{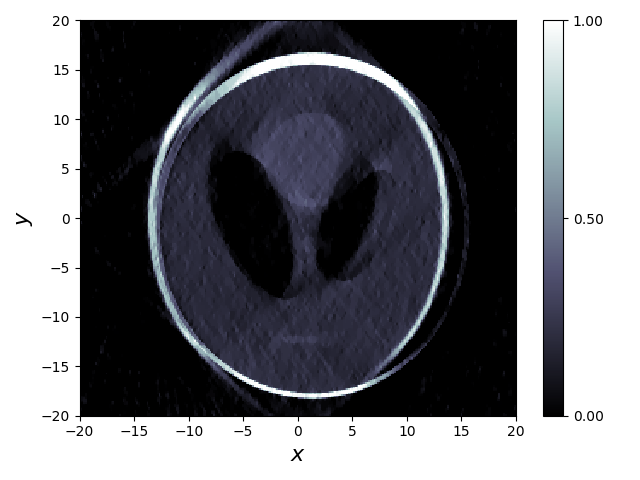}}%
  \hfil
  \subfloat[\label{fig:L2_c}]{\includegraphics[width=0.33\textwidth]{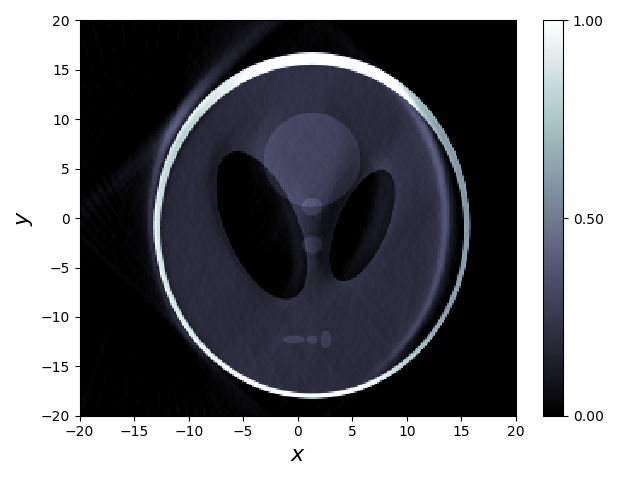}}%
  \caption{Reconstructions using $\ell_2$ prior with different regularization parameters:\newline (a) $\gamma=1$, (b) $\gamma=100$, and (c) $\gamma=10\,000$.}%
  \label{fig:results_SL_L2}%
\end{figure}

\begin{remark}\label{rem:big}
These inverse problems are highly underdetermined and ill-posed. The total number of pixels is \imDim$^2 =$ \imDimTot\, but the number of data points in the examples are only \parLines$\cdot$\Angles\,$ = $\,\nrMeasurements\ and \parLines$\cdot$\AnglesHand\,$ = $\,\nrMeasurementsHand, respectively. Also note that solving the corresponding optimal transport problems explicitly  would amount to matrices of sizes \imDim$^2 \times$\imDim$^2$, which means solving linear programs with over \nrVariables\, variables. 
\end{remark}

\begin{figure}%
  \centering
  \hfil
  \subfloat[\label{fig:phantom2fbp}]{\includegraphics[width=0.5\textwidth]{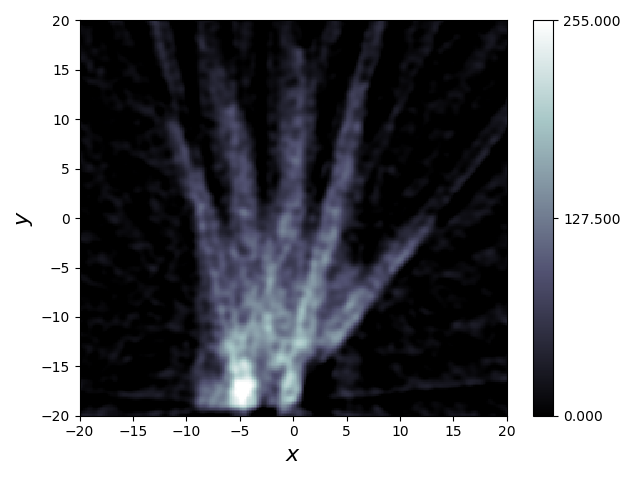}}%
  \hfil
  \subfloat[\label{fig:TVHand}]{\includegraphics[width=0.5\textwidth]{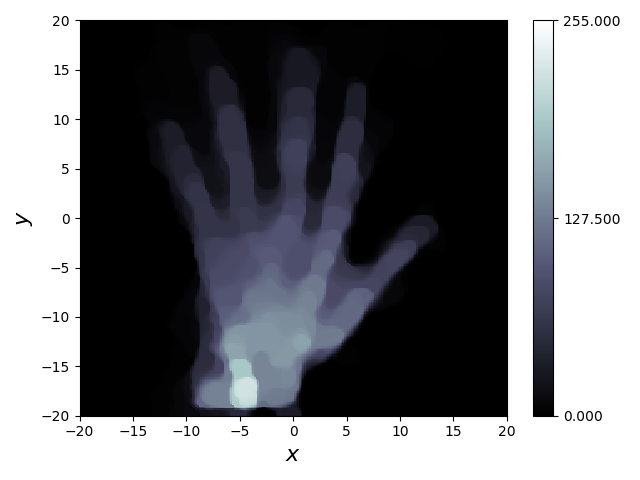}}%
  \hfil
  \subfloat[\label{fig:L2Hand}]{\includegraphics[width=0.5\textwidth]{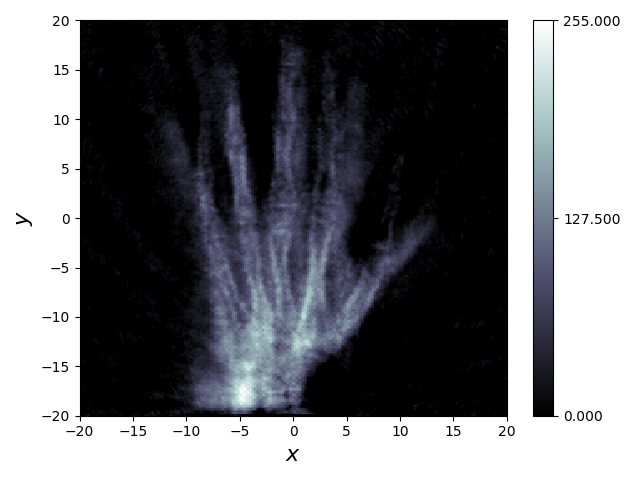}}%
  \hfil
  \subfloat[\label{fig:OMTHand}]{\includegraphics[width=0.5\textwidth]{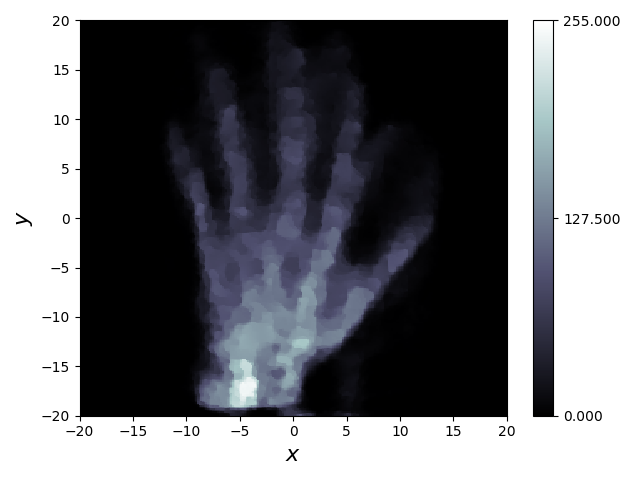}}%
  \caption{Reconstructions using different methods. (a) Filtered backprojection, (b) reconstruction using TV-regularization, (c) reconstruction with $\ell_2^2$-prior and TV-regularization, and (d) reconstruction with optimal transport prior and TV-regularization.}%
  \label{fig:resultsHand}%
\end{figure}

\subsection{Reconstruction of Shepp-Logan image}
The reconstructions are shown in Figure~\ref{fig:results}.
Both the the filtered packprojection reconstruction in Figure~\ref{fig:fbp} and the 
TV-reconstruction in Figure~\ref{fig:TV} suffers from artifacts and sever vertical blurring due to poor vertical resolution resulting from the limited angle measurements. 
Figure \ref{fig:L2} shows the reconstruction with $\ell_2$-prior. Some details are visible, however these are at the same locations as in the prior and does not adjust according to the measurements of the phantom. Considerable artifacts also appear in this reconstruction, typically as fade-in-fade-out effects where the prior and the data do not match. 
The fade-in-fade-out effect that often occur when using strong metrics for regularization is illustrated in Figure~\ref{fig:results_SL_L2}. Selecting a low value $\gamma$ (Figure~\ref{fig:L2_a}) results in a reconstruction close to the TV-reconstruction and selecting a large value $\gamma$ gives a reconstruction close to the prior (Figure~\ref{fig:L2_c}).
By selecting a medium value $\gamma$ one gets a reconstruction that preserves many of the details found in the prior, however they remain at the same position as in the prior and hence are not adjust to account for the measurement information.

The reconstruction with optimal mass transport prior is shown in Figure~\ref{fig:OMT}. Some blurring occurs, especially in the top and the bottom of the image, however the overall shape is better preserved compared to the other reconstructions. Fine details are not visible, but the major features are better estimated compared to the TV- and $\ell_2$-reconstructions. 
This example illustrates how one can improve the reconstruction by incorporating prior information, but without the fade-in-fade-out effects that typically occurs when using a strong metric such as $\ell_2$ for regularization.

\begin{figure}%
  \centering
  \hfil
  \subfloat[\label{fig:omt_masked_prior}]{\includegraphics[width=0.5\textwidth]{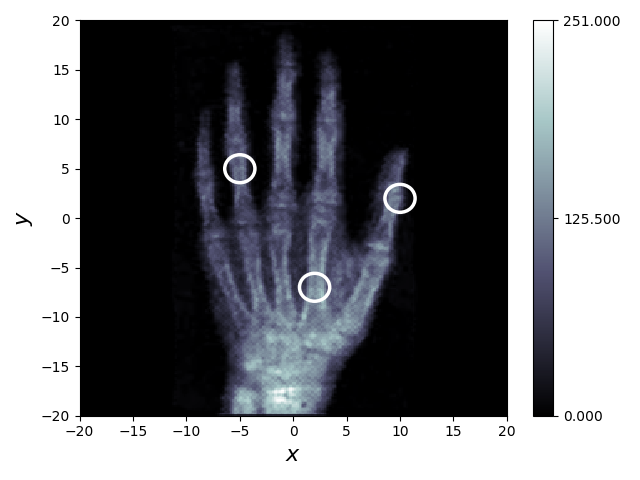}}%
  \hfil
  \subfloat[\label{fig:omt_mass_trans}]{\includegraphics[width=0.5\textwidth]{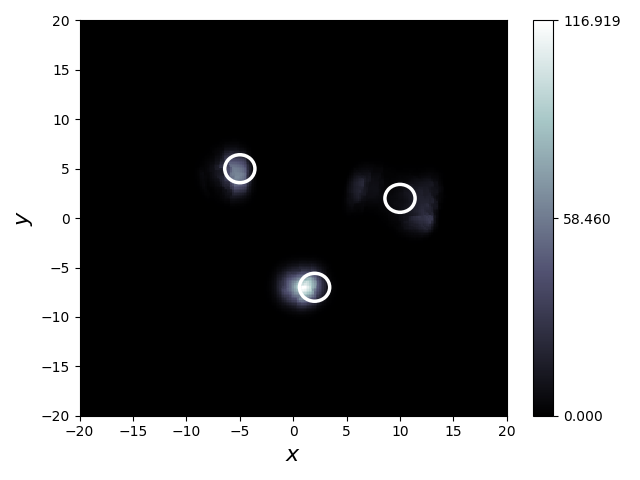}}%
  \hfil
  \subfloat[\label{fig:omt_recon}]{\includegraphics[width=0.5\textwidth]{\myFigPathHand \HandOMT.png}}%
  \caption{Visualizing how the mass from certain regions are transported. (a) the prior and three regions that are considered, (b) how the mass from certain regions in the prior are transported to the reconstruction, and (c) the reconstruction using optimal mass transport prior.}%
  \label{fig:results_mass_transport}%
\end{figure}

\subsection{Reconstruction of hand image}

The reconstructions are shown in Figure~\ref{fig:resultsHand}, and 
we obtain similar results as in the previous example.
The filtered backprojection reconstruction in Figure~\ref{fig:phantom2fbp} is quite fragmented, and although some details are visible there is also a considerable amount of noise and artifacts
 throughout the image. The TV-reconstruction, shown in Figure~\ref{fig:TVHand}, is smeared and few details are visible. Moreover, in the reconstruction with $\ell_2$ prior, shown in Figure~\ref{fig:L2Hand}, similar artifacts as in the Shepp-Logan example are present. Again, they steam from a miss-match between the prior and the data. Note especially that the thumb, but also the index and middle finger, almost look dislocated or broken in the reconstruction. 
The reconstruction with optimal mass transport prior is shown in Figure~\ref{fig:OMTHand}. Also in this case the fine details are not visible. Note however that details are more visible in \ref{fig:OMTHand}, compared to the TV-reconstruction \ref{fig:TVHand}, and that the reconstruction in \ref{fig:OMTHand} does not suffer from the the same kind of artifacts as the $\ell_2$-regularized reconstruction in \ref{fig:L2Hand}. 

To illustrate the effect of the optimal mass transport prior on the final reconstruction, we also include Figure~\ref{fig:results_mass_transport}. 
Figure~\ref{fig:omt_masked_prior} shows the prior \ref{fig:priorHand} with certain areas marked, and  Figure~\ref{fig:omt_mass_trans} shows how the mass from these areas in \ref{fig:omt_masked_prior} are transported to \ref{fig:omt_recon}.
For reference the  optimal mass transport reconstruction from \ref{fig:OMTHand} is shown in Figure~\ref{fig:omt_recon}. By comparing the images one can see that the thumb in the prior is to a large extent transported to the location of the thumb in the phantom, which is what the optimal mass transport is intended to do. However, a fraction of the mass from the thumb is also transported to the index finger, giving rise to some artifacts.
One can also note that there is a certain ``halo-effect'' around each region, especially  for the region on the thumb. This is most likely due to the computational regularization $ \epsilon D(M)$ in \eqref{eq:Teps}, which forces the optimal $M$ to be strictly positive. It would therefore be desirable to precondition or in other ways improve the numerical conditioning of the steps in Algorithm~\ref{alg:prox_Sinkhorn}, to allow for the computation of the proximal operator for lower values of $\epsilon$.

\begin{remark} \label{rem:otherTom}
As mentioned in the beginning of this section, the papers \cite{abraham2017tomographic, benamou2015iterative} also consider tomography problems with optimal mass transport priors. In both these papers, the optimization problems are on the form
\[
\min_{\mue} T(\mu_0, \mue)+ \sum_{\ell=1}^L \gamma_\ell\tilde T(P_{\theta_\ell}\mue, b_\ell),
\]
where $T$ and $\tilde T$ are transportation costs and $P_{\theta_\ell}$ denotes the projection along the angle $\theta_\ell$ (i.e., a partial Radon transform).
In our formulation we use a hard data matching constraint based on the $L_2$-norm (see \eqref{eq:OMT}) which models Gaussian noise. A computational approach is provided in \cite{benamou2015iterative} that also builds on Sinkhorn iterations (Bregman projections). However, the computational procedure builds on that the adjoint $P_{\theta_\ell}^T$ commutes with elementwise functions, which is a result of the implementation of the partial ray transform $P_{\theta_\ell}$ using nearest neighbor interpolation.%
\footnote{Implementing the partial ray transform $P_{\theta_\ell}$ using nearest neighbor interpolation results in a matrix representation where each column is an elementary unit vector.}  On the other hand, our method allows for an arbitrary (linear) forward operator. 
\end{remark}

\section{Concluding remarks and further directions}\label{sec:conclusions}
In this work we have considered computational methods for solving inverse problems containing entropy regularized optimal mass transport terms. First, using a dual framework we have generalized the Sinkhorn iteration, thereby making it applicable to a set of optimization problems. In particular, the corresponding proximal operator of the entropy regularized optimal mass transport cost is computed efficiently. Next, we use this to address a large class of inverse problems using variable splitting. In particular 
we use a Douglas-Rachford type method to solve two problems in computerized tomography where prior information is incorporated using optimal mass transport.

Interestingly, both the Sinkhorn iterations and the proposed approach for computing the proximal operator are identical to coordinate ascent of the dual problem. For these problems the coordinate ascent step can be computed explicitly by \eqref{eq:sinkhorn} or \eqref{eq:equations_prox}, respectively. In this setting the hard constraint \eqref{eq:T_dualb} is replaced by the soft constraint \eqref{eq:barrier}. The Hessian of the barrier term is given explicitly by 
\[
-\frac{1}{\epsilon}\begin{pmatrix}
\diag(u_0\odot (Ku_1)) & \diag(u_0) K\diag(u_1)\\\diag(u_1) K^T\diag(u_0)&\diag(u_1\odot (K^Tu_0))
\end{pmatrix}
\]
and multiplication of this can be computed quickly (cf. Remark~\ref{rem:Fast_K}). This opens up for accelerating convergence  using, e.g., quasi-Newton or parallel tangent methods \cite{luenberger1984linear}.  

The methodology presented in this paper naturally extends to problems with several optimal mass transport terms.
In particular this would be useful for cases where estimates are not guaranteed to have the same mass or where some parts are not positive. An examples of such a problem is the computation of a centroid (barycenter) from noisy measurements (cf. \cite{chizat2016scaling}), i.e.,  
\begin{equation}
\min_{\mu_{\ell},\,\ell=0,\ldots, L} \;\sum_{\ell=1}^L \left(T_\epsilon(\mu_0,\mu_\ell)+ \frac{1}{2\sigma_\ell}\|\mu_\ell-\nu_\ell\|_2^2\right).
\end{equation}
This problem has applications in clustering and will be considered in future work.

\section{Acknowledgements}
The authors would like to thank Yongxin Chen and Tryphon Georgiou  for input and fruitful discussions on the optimal mass transport problem, and Krister Svanberg for insightful discussion on optimization. The authors would also like to thank Jonas Adler and Ozan \"Oktem for input on the example and the implementation in ODL. The figures in Figure~\ref{fig:phantomHand} and Figure~\ref{fig:priorHand} have generously been provided by Klas Modin, and can be found in \cite{bauer2015diffeomorphic}. We also would like to thank the anonymous referees valuable feedback and for bringing the references \cite{peyre2015entropic, cuturi2015asmoothed, schmitzer2016stabilized} to our attention.

\appendix

\section{Proof of Proposition \ref{prp:dual}} \label{app:A}
We seek the dual problem of the following optimization problem
\begin{align}
\min_{M\ge 0, \mue} \qquad & \tr(C^TM)+\epsilon D(M)+g(\mue)\label{eq:Pgen_eps}\\
\mbox{subject to } \quad&  \mu_0=M \ett_{n_1}\nonumber\\
& \mue=M^T \ett_{n_0}.\nonumber
\end{align}
Lagrange relaxation gives the Lagrangian
\begin{align*}
L(M,\mue, \lambda_0, \lambda_1)&=\tr(C^TM)+\epsilon D(M)+g(\mue)\\&\quad+\lambda_0^T(\mu_0-M \ett_{n_1})
+\lambda_1^T(\mue-M^T \ett_{n_0})\\
&=\tr((C-\ett_{n_1}\lambda_0^T-\lambda_1\ett_{n_0}^T)^TM)+\epsilon D(M)\\&\quad+\lambda_0^T\mu_0
-((-\lambda_1)^T\mue-g(\mue)). 
\end{align*}
First note that by definition $g^*(-\lambda_1)=\max_{\mue}(-\lambda_1)^T\mue-g(\mue)$, hence the last term equal $-g^*(-\lambda_1)$ for the minimizing $\mue$. Next, the minimizing $M$ is unique and satisfies
\[
\epsilon \log (m_{ij})+c_{ij}-\lambda_0(i)-\lambda_1(j)=0
\]
or equivalently 
\[
M= \diag(\exp(\lambda_0/\epsilon))\exp(-C/\epsilon)\diag(\exp(\lambda_1/\epsilon)).
\]
By plugging these into the Lagrangian we obtain
\begin{align*}
&\min_{M, \mue} L(M,\mue, \lambda_0, \lambda_1)= \epsilon\sum_{i=1}^{n_0}\sum_{j=1}^{n_1}(1-m_{ij})+\lambda_0^T \mu_0-g^*(-\lambda_1)\\
&= \lambda_0^T\mu_0-g^*(-\lambda_1)-\epsilon \exp(\lambda_0^T/\epsilon)\exp(-C/\epsilon)\exp(\lambda_1/\epsilon) +\epsilon n_0n_1, 
\end{align*}
which is the objective function of the dual problem. 

Moreover, by Assumption~\ref{ass:g} and the fact that $T_\epsilon(\mu_0, \mue)<\infty$ restricts $\mue$ to a compact set, there exists an optimal solution to \eqref{eq:Pgen_eps}. 
To see this, note that if $\mue$ is a point such that $g(\mue) < \infty$ and $\sum_{i=1}^{n_0} \mu_0(i) = \sum_{j=1}^{n_1} \mue(j)$, then $M = \mu_0 \mue^T$ is an inner point where the objective function takes a finite value. 
From this it also follows that the minimum of \eqref{eq:Pgen_eps} is finite. 
Therefore by \cite[Proposition 5.3.2]{bertsekas1999nonlinear} strong duality holds.$\hfill\square$

\section{Proof of Theorem \ref{thm:main}} \label{app:B}
First we would like to show that Algorithm~\ref{alg:prox_Sinkhorn} corresponds to coordinate ascent of the dual problem, i.e., that the steps \ref{alg:step_sinkhorn} and \ref{alg:step_prox} corresponds to the maximization of $\lambda_0$ and $\lambda_1$, respectively. This follows directly  for step \ref{alg:step_sinkhorn} since it is identical to \eqref{eq:dual_stat2a} in Lemma \ref{lm:grad_cond}. Note that step \ref{alg:step_sinkhorn} is the Sinkhorn iterate \eqref{eq:sinkhorn_a} (cf. Corollary~\ref{cor:Sinkhorn}).

Next, consider the condition \eqref{eq:dual_stat2b} in Lemma \ref{lm:grad_cond}, from which it follows that the minimizing $\lambda_1$ satisfies
\begin{align*}
\mu_1-\sigma\lambda_1&=\exp(\lambda_1/\epsilon)  \odot  \left(K^T u_0 \right),
\end{align*}
where we let $K=\exp(-C/\epsilon)$ and $u_0=\exp(\lambda_0/\epsilon)$. Taking the logarithm and adding $(\mu_1-\sigma\lambda_1)/(\sigma\epsilon)-\log(\sigma\epsilon)$ to each side, we get
\begin{align*}
\frac{\mu_1-\sigma\lambda_1}{\sigma\epsilon}+\log\left(\frac{\mu_1-\sigma\lambda_1}{\sigma\epsilon}\right)&=\frac{\mu_1}{\sigma\epsilon}  +\log\left(\frac{K^T u_0}{\sigma\epsilon} \right),
\end{align*}
or equivalently
\begin{align}\label{eq:wright_omega_exp}
\frac{\mu_1-\sigma\lambda_1}{\sigma\epsilon}=\omega\left(\frac{\mu_1}{\sigma\epsilon}  +\log\left(\frac{K^T u_0}{\sigma\epsilon} \right)\right),
\end{align}
where $\omega$ denotes the elementwise Wright omega function, i.e., the function mapping $x\in\mR$ to $\omega(x)\in \mR_+$ for which $x=\log(\omega(x))+\omega(x)$ \cite{corless2002wright}. The function is well defined as a function $\mR\to\mR_+$ which is the domain and range of interest in our case. This expression of $\lambda_1$ is equivalent to step \ref{alg:step_prox} (and \eqref{eq:dual_prox2}), and we have thus shown that Algorithm~\ref{alg:prox_Sinkhorn} is a coordinate ascent algorithm for the dual problem \eqref{eq:dual_prox}.

Next note that \eqref{eq:dual_prox} is a strictly convex optimization problem. This follows since the Hessian 
\begin{equation}\label{eq:Hessean_prox}
-\frac{1}{\epsilon}\begin{pmatrix}
\diag(u_0\odot (Ku_1)) & \diag(u_0) K\diag(u_1)\\ \diag(u_1) K^T\diag(u_0)&\diag(u_1\odot (K^Tu_0))
\end{pmatrix}-\begin{pmatrix}
0 & 0\\ 0&\sigma I
\end{pmatrix}
\end{equation}
is strictly negative definite. This can be seen by noting that the first term is diagonally dominant, hence negative semidefinite. Since $\sigma>0$, then any zero eigenvector of \eqref{eq:Hessean_prox} can only have non-zero elements in the first block. However, the $(1,1)$-block of \eqref{eq:Hessean_prox}, i.e., $-\diag(u_0\odot (Ku_1))/\epsilon$, is negative definite, hence no zero eigenvalue exists and \eqref{eq:Hessean_prox} is strictly negative definite. 

As seen above, the optimization problem is strictly convex, hence there is a unique stationary point which is also the unique maximum.
To show convergence, note that the objective function in \eqref{eq:dual_prox} is continuously differentiable, hence any limit point of the coordinate ascent is stationary  \cite[Proposition 2.7.1]{bertsekas1999nonlinear}. 
Further, since the suplevel sets of \eqref{eq:dual_prox} are bounded and there is a unique stationary point, Algorithm~\ref{alg:prox_Sinkhorn} converges to unique maximum. 
Finally, locally linear convergence follows from \cite[Theorem 2.2]{bezdek1987local} since the optimization problem is strictly convex. $\hfill\square$\\

\begin{remark}
For future reference, we note that \eqref{eq:wright_omega_exp} can be rewritten as
\begin{align}\label{eq:wright_omega_exp_2}
\frac{\lambda_1}{\epsilon}=&\;
\frac{\mu_1}{\sigma\epsilon}-\omega\left(\frac{\mu_1}{\sigma\epsilon}  +\log\left(\frac{K^T u_0}{\sigma\epsilon} \right)\right)\nonumber\\
=&\;\frac{\mu_1}{\sigma\epsilon}-\left(\frac{\mu_1}{\sigma\epsilon}  +\log\left(\frac{K^T u_0}{\sigma\epsilon}\right)-\log\left[\omega\left(\frac{\mu_1}{\sigma\epsilon}  +\log\left(\frac{K^T u_0}{\sigma\epsilon} \right)\right)\right]\right)\\
=&\;\log\left(
\frac{
\sigma\epsilon\omega\left(\frac{\mu_1}{\sigma\epsilon}  +\log\left(\frac{K^T u_0}{\sigma\epsilon} \right)\right)
}{
K^T u_0
}
\right),\nonumber
\end{align}
where we in the second equality use the definition $\omega(x)=x-\log(\omega(x))$. This expression can alternatively be used in Algorithm \ref{alg:prox_Sinkhorn}. However, our experience is that \eqref{eq:wright_omega_exp} is better conditioned than \eqref{eq:wright_omega_exp_2}.
\end{remark}

\section{Connection with method based on Dykstra's algorithm} \label{app:C}
We would like to thank the reviewers for pointing out the reference \cite{peyre2015entropic}. In fact Algorithm \ref{alg:prox_Sinkhorn} can be seen as a (nontrivial) special case of the iterations in \cite[Proposition 3.3]{peyre2015entropic}. We will here provide a separate derivation leading to a simplified but equivalent version of this algorithm. 
The algorithm \cite[Proposition 3.3]{peyre2015entropic} also address the optimization problem \eqref{eq:Pgen_eps0} and can be used when the entropic proximal, defined by 
\begin{equation*}
\Prox^{{\rm \overline{KL}}}_{\sigma g}(z):=\argmin_{x}\; \sigma g(x)+D(x|z),
\end{equation*}
where $D(x|z)=\sum_{i=1}^{n} (x_{i}\log(x_{i}/z_{i})-x_{i}+z_i)$,
is fast to compute. It can be noted that computing the entropic proximal is equivalent to
\begin{align*}
\Prox^{{\rm \overline{KL}}}_{\sigma g}(z):= & \; \argmin_{x,x'} \qquad g(x)+\frac{1}{\sigma}D(x'|z)\\
& \;\; \text{subject to} \quad x=x',
\end{align*}
which has a Lagrange dual that is given by (cf. the proof of Proposition \ref{prp:dual})
\[
\max_{\lambda}\; -g^*(-\lambda)-\exp(\sigma\lambda^T)z/\sigma+\ett^T z/\sigma,
\]
and where the optima of the primal and dual problems relate as $\sigma\lambda=\log(x./z)$.
The maximizing argument $\lambda$ is specified by
$$0\in \partial g^* (-\lambda)-\exp(\sigma\lambda)\odot z/\sigma.$$
Noting that this condition is equivalent to \eqref{eq:dual_stat2b} with $z/\sigma= \exp(-C^T/\epsilon)\exp(\lambda_0/\epsilon)$ and $\epsilon=1/\sigma$, the dual update of $\lambda_1$ can be written as
\begin{equation}\label{eq:peyre}
\exp(\lambda_1/\epsilon) = \frac{\Prox^{{\rm \overline{KL}}}_{\epsilon^{-1} g}(\exp(-C^T/\epsilon)\exp(\lambda_0/\epsilon)) }{\exp(-C^T/\epsilon)\exp(\lambda_0/\epsilon)}=\frac{\Prox^{{\rm \overline{KL}}}_{\epsilon^{-1} g}(K^Tu_0) }{K^Tu_0},
\end{equation}
where as before $K=\exp(-C/\epsilon)$ and $u_i=\exp(\lambda_i/\epsilon)$, for $i=0,1$.
Using this we can state the simplified, but equivalent, version of \cite[Proposition 3.3]{peyre2015entropic}, shown in Algorithm \ref{alg:prox_Peyre}. The algorithm is equivalent in the sense that with $K=\xi^T$, then it holds that $b^{(2\ell)}=b^{(2\ell-1)}=u_0^{(\ell)}$ and $a^{(2\ell+1)}=a^{(2\ell)}=u_1^{(\ell)}$ for $\ell\ge 1$. Similar algorithms have also been derived in the recent preprints 
\cite{chizat2016scaling, schmitzer2016stabilized}.

\begin{algorithm}
\caption{Simplified version of \cite[Proposition 3.3]{peyre2015entropic}.}
\label{alg:prox_Peyre}
\begin{algorithmic}[1]
\REQUIRE $\epsilon$, $C$, $\lambda_0$, $\mu_0$, $g$
\STATE $K=\exp(-C/\epsilon)$ and $u_1=\ett$
\STATE $\ell=0$
\WHILE{Not converged}
\STATE $\ell \leftarrow \ell+1$
\STATE $u_0^{(\ell)} \leftarrow \mu_0./(Ku_1^{(\ell-1)})$ 
\STATE $\displaystyle u_1^{(\ell)} \leftarrow \frac{\Prox^{{\rm \overline{KL}}}_{\epsilon^{-1} g}(K^Tu_0^{(\ell)}) }{K^Tu_0^{(\ell)}}$
\ENDWHILE
\ENSURE $\mue \leftarrow  u_1^{(\ell)} \odot (K^T u_0^{(\ell)})$
\end{algorithmic}
\end{algorithm}

The case which is a main focus for this paper is $g(\cdot)=\frac{1}{2\sigma}\|\cdot - \mu_1\|^2_2$, and comparing the expressions \eqref{eq:peyre} and \eqref{eq:wright_omega_exp_2} indicates that
\[
\Prox^{{\rm \overline{KL}}}_{\epsilon^{-1} g}(K^Tu_0)=\sigma\epsilon\,\omega\left(\frac{\mu_1}{\sigma\epsilon}  +\log\left(\frac{K^T u_0}{\sigma\epsilon} \right)\right).
\]
This can in fact be verified by direct computations along the lines of the proof of Theorem~\ref{thm:main}.

\section{Parameters in the numerical examples} \label{app:D}

The problem is set up with noise level $5\%$ in the Shepp-Logan example and $3\%$ in the hand example. The parameter $\kappa$ is selected to be $120\%$ of the norm of the noise. That is, 
in the Shepp logan example white noise is generated and normalized so that  $\|b-A\mut\|_2/\|A\mut\|_2=0.05$ and $\kappa=1.2 \cdot 0.05\cdot\|A\mut\|_2$. This ensures that the true image $\mut$ belong to the feasible region $\{\mu:\,\|b-A\mu\|_2\le \kappa\}$.
The cost function in the optimal mass transport distance is given by
\begin{equation}\label{eq:costfucntion}
c(x_1,x_2)=\min(\|x_1-x_2\|_2, 20)^2,
\end{equation}
where the truncation at $20$ is done in order to improve the conditioning of the computations. Further, the proximal operator of the optimal mass transport functional is computed using \SinkhornIter\, generalized Sinkhorn iterations.
Each optimization problem is solved using \DRiter\, iterations in the Douglas-Rachford algorithm. The step size parameters $\sigma_i$ are set to be 
$\sigma_i=(\tau \|L_i\|_{\rm op}^2)^{-1}$
where $\|L\|_{\rm op}=\sup_{\|x\|_2\le 1} \|L x\|_2$ is the operator norm (approximated using the power iterations in ODL \cite{adler2017ODL}). The remaining parameters are selected according to Table~\ref{tab:param_val}. The gradient operator used in the TV-terms is the default \texttt{Gradient} operator in ODL, which pads the boundaries with zeros and applies a forward-difference in each dimension (see the documentation \url{http://odlgroup.github.io/odl/}). Other options are available, e.g., zero-order-hold padding 
or periodic boundary conditions, as well as backward- or central-difference. In our case, zero-padding is not an issue since both phantoms are zero along the boundary, see Figure~\ref{fig:phantom} and Figure~\ref{fig:phantomHand}. 
For the filtered backprojection reconstruction we use the ODL-implementation with a Hann filter with filter parameter $0.7$ for the Shepp-Logan example and $0.5$ for the hand example.

\noindent 
\begin{table}[ht]
\caption{\label{tab:param_val} Parameter values for the variational problems and reconstruction algorithms. A $\star$  means that the parameter is not used in this problem.}
\centering
\begin{tabular}{l|l| cc |cc}
\toprule
\multicolumn{2}{c|}{Reconstruction example} & \multicolumn{2}{c|}{Parameters in} & \multicolumn{2}{c}{Parameters in} \\
\multicolumn{2}{c|}{} & \multicolumn{2}{c|}{optimization problem} & \multicolumn{2}{c}{Algorithm \ref{alg:DR}} \\
\midrule
Phantom     & Objective function & $\gamma$                    & $\epsilon$ & $\tau$ & $\lambda$ \\
\midrule
Shepp-Logan & TV             & $\star$                     & $\star$    & $0.05$ & $1$ \\
            &  TV + $L_2^2$   & $\{1, 10, 100, 10 \, 000\}$ & $\star$    & $0.05$ & $1$  \\
            & TV + OMT       & $4$                         & $1$        & $5$    & $1.8$ \\
\midrule
Hand        & TV             & $\star$                     & $\star$    & $1/\sqrt{2}$  & $1$ \\
            & TV + $L_2^2$     & $10$                        & $\star$    & $1/\sqrt{2}$  & $1$  \\
            &  TV + OMT        & $4$                         & $1.5$      & $500\sqrt{2}$  & $1.8$ \\
\bottomrule
\end{tabular}
\end{table} 

\bibliographystyle{plain}

\bibliography{bib_johan}

\end{document}